\documentclass[10pt]{article}
\usepackage{amsmath,amssymb,amsthm,mathrsfs,calc,color}
\usepackage{graphicx}
\usepackage[british]{babel}
\usepackage{enumerate}

\newcommand{\PPP}[1]{\operatorname{\mathbb{P}}\left(\,#1\,\right)}
\newcommand{\EEE}[1]{\operatorname{\mathbb{E}}\left[\,#1\,\right]}

\newtheorem {theorem}{Theorem}[section]
\newtheorem {proposition}[theorem]{Proposition}
\newtheorem {lemma}[theorem]{Lemma}
\newtheorem {corollary}[theorem]{Corollary}

{\theoremstyle{definition}

}
{\theoremstyle{theorem}

}
\newcommand{\conv}[2][n]{\underset{#1\rightarrow #2}{\longrightarrow}}

\newcommand{\ind}[1]{\mathbf{1}_{#1}\,}

\def\NN{\mathbb{N}}
\def\PP{\mathbb{P}}
\def\QQ{\mathbb{Q}}
\def\RR{\mathbb{R}}

\usepackage[usenames,dvipsnames]{xcolor}

\newcommand{\hide}[1]{}

\usepackage[normalem]{ulem}
\numberwithin{equation}{section}
\usepackage[french]{varioref}

\newenvironment{prooft}[1]{\noindent {\bf Proof of #1}} {\hfill $\square$  \noindent}

\begin{document}

\title{The largest order statistics for the inradius in an isotropic STIT tessellation}
\author{Nicolas Chenavier \footnote{Universit\'e du Littoral C\^ote d'Opale, LMPA Joseph Liouville, BP 699, F-62228 Calais Cedex, France; Email address: nicolas.chenavier@univ-littoral.fr} and Werner Nagel\footnote{Friedrich-Schiller-Universit\"at Jena, Fakult\"at f\"ur Mathematik und Informatik, D-07737 Jena, Germany. Email address: werner.nagel@uni-jena.de } }

\maketitle
\date{}

\begin{abstract}
A planar stationary and isotropic STIT tessellation at time $t>0$ is observed in the window $W_\rho={t^{-1}}\sqrt{\pi \ \rho}\cdot [-\frac{1}{2},\frac{1}{2}]^2$,  for $\rho>0$. With each cell of the tessellation, we associate the inradius, which is the radius of the largest disk contained in the cell. Using the Chen-Stein method, we compute the limit distributions of the largest order statistics for the inradii of all cells whose nuclei are contained in $W_\rho$ as $\rho$ goes to infinity.
\end{abstract}

\noindent\textbf{Keywords:} Stochastic Geometry; Random Tessellations; Extreme Values; Poisson Approximation


\noindent \textbf{AMS 2010 Subject Classifications:} 60D05 . 60G70 . 60F05 . 62G32

\section{Introduction}
{A planar tessellation $T$ in $\mathbb{R}^{2}$ is a countable collection of non-empty convex
polygons, called \textit{cells}, with disjoint interiors which subdivides the plane $\RR^2$, and such that the number of cells intersecting any
bounded subset of $\mathbb{R}^{2}$ is finite.} Let  $\mathcal{T}$ be the set of {planar} 
tessellations, endowed with the usual $\sigma $-field, see {p. \pageref{page:tessellation}}. A random {planar} tessellation  is a random variable with values in $\mathcal{T}$. {For a complete account on random planar tessellations, and more generally random tessellations in a $d$-dimensional Euclidean space, $d\geq 2$}, we refer to the books  \cite{SW, SKM}.

One of the important models is the STIT (STable under ITerations) tessellation. This random tessellation was introduced in \cite{NagelWeiss05} and has potential applications for the modeling of crack patterns or of fracture structures, such as the so-called craquel\'ee on pottery surfaces, drying films of colloidal particles, or structures observed in geology \cite{LP,MoMa}. There are quite a few theoretical results for STIT tessellations, including  ergodicity \cite{MN3}, mixing properties \cite{lachiezerey}, computations of distributions \cite{NagelWeiss11} and  a Mecke-type formula \cite{NNTW17}.

Let ${\bf Y}=(Y_t,\, t>0)$ be a STIT tessellation process in the Euclidean plane $\RR^2$, where for any $t>0$ the tessellation $Y_t$ is spatially stationary and isotropic. The Euclidean plane $\RR^2$ is endowed with its Euclidean norm $||\cdot ||$.   Let $t>0$ be fixed.
 With each cell $z\in Y_t$ we associate the incenter $c(z)$ of $z$, defined as the center of the largest disk included in $z$. For each cell the incenter exists and is almost surely unique. Now, let $B$ be a Borel subset in $\RR^2$ with area $a(B)\in (0,\infty)$. The \textit{cell intensity} $\gamma_t$  of $Y_t$ is defined as the mean number of cells per unit area, i.e.
 \[\gamma_t=\frac{1}{a(B)}\EEE{\#\{z\in Y_t: c(z)\in B\}}.\] 
Because a STIT tessellation is stationary, the intensity is in fact independent of the Borel subset $B\subset\RR^2$. In the present paper, the STIT tessellation $Y_t$ is scaled in such a way that $\gamma_t=\frac{t^2}{\pi}$. 
The \textit{typical cell} is a random polygon $\mathcal{Z}$ whose distribution is given by 
\begin{equation} \label{eq:deftypical} 
\EEE{g(\mathcal{Z})} = \frac{1}{\gamma_t \, a(B)}\EEE{\sum_{z\in Y_t: c(z)\in B}g(z-c(z))},
\end{equation}
for all  nonnegative measurable functions $g:\mathcal{K}_0 \rightarrow\RR$ on the set $\mathcal{K}_0$ of centered convex bodies, i.e. all non-empty convex compact sets  $K\subset \RR^2$ with $c(K)=0$, endowed with the Hausdorff topology.
(For those random tessellations for which the incenters of cells are not a.s. unique, choose another center function $c$ for the cells).

It is notable  that the distribution of the typical cell in a STIT tessellation is the same as in a stationary Poisson line tessellation with corresponding parameters \cite{NagelWeiss05}, i.e. in our case, the stationary and isotropic Poisson line tessellation with the same cell number intensity $\gamma_t$. 

In the present paper, as a first approach on extremal properties of a STIT tessellation, we study the largest order statistics for the inradius. The inradius  is one of the rare geometric characteristics for which the distribution can be made explicit. Before stating our main theorem, we introduce some notation. For each cell $z\in Y_t$, we denote by  $R(z)$ the \textit{inradius} of $z$, i.e. the radius of the largest disk included in the cell $z$. It is known that the inradius of the typical cell of $Y_t$ has an exponential distribution with parameter $2t$ (Lemma 3 in \cite{NagelWeiss05}), i.e.
\begin{equation}\label{eq:exponential} 
\PPP{R(\mathcal{Z})>v}=e^{-2tv}, \quad v\geq 0.
\end{equation}
Consider the family $W_\rho=t^{-1} \sqrt{\pi \ \rho}\cdot [-\frac{1}{2},\frac{1}{2}]^2$, $\rho>0$,  of squares which we refer to as windows. Now, for each threshold $v\geq 0$, let $N_{W_\rho}(v)$ be the number of cells with incenters in $W_\rho$ and with inradii larger than $v$, i.e.
\begin{equation}\label{eq:defNWvrho}
N_{W_\rho}(v):=\sum_{z\in Y_t: c(z)\in W_\rho}\ind{R(z)>v}.
\end{equation}
The cells with inradius exceeding $v$ are referred to as \textit{exceedances}. According to \eqref{eq:deftypical}, the mean number of exceedances is given by $\EEE{N_{W_\rho}(v)} = \gamma_t \,  t^{-2}\,   \pi \, \rho \  \PPP{R(\mathcal{Z})>v}$.

In this paper, we provide a Poisson approximation of the number of exceedances when the size $\rho$ of the window  goes to infinity. The underlying threshold which we consider has to depend on $\rho$. To define it, let $\tau>0$ be a fixed value. The threshold is chosen as a function $v_\rho$, $\rho>0$, in such a way that the mean number of exceedances equals $\tau$, i.e.
\begin{equation} \label{Campbellmark}
\EEE{N_{W_\rho}(v_\rho)} = \gamma_t \,  t^{-2}\, \pi \,\rho \,  \PPP{R(\mathcal{Z})>v_\rho}=\tau.\end{equation}
Thus, the threshold is chosen as:
 \begin{equation}\label{eq:defvrho}
v_\rho:=v_\rho(\tau)=\frac{1}{2t}(\log\rho-\log\tau).
\end{equation} 
We consider the convergence of distributions with respect to the total variation. In the particular case of two random variables $X$ and $Z$  with values in $\NN_0 =\{ 0,1,2,\ldots \}$, we recall that the total variation distance is given by
\[d_{TV}(X,Z) := 2\sup_{A\subset\NN_0}|\PPP{X\in A} - \PPP{Z\in A}|.\] We are now prepared to state our main result.
\begin{theorem}
\label{Th:extremes}
Let $Y_t$ be a stationary and isotropic STIT tessellation in $\RR^2$ at time $t>0$ with cell number intensity $\gamma_t=\frac{t^2}{\pi}$. Let $\tau>0$ be fixed. Let $v_\rho$ be as in \eqref{eq:defvrho} and $N_{W_\rho}(v_\rho)$ as in \eqref{eq:defNWvrho}. Furthermore, let $Z$ be a random variable with a Poisson distribution with parameter $\tau$. Then
\begin{equation}\label{eq:convTV}
d_{TV}(N_{W_\rho}(v_\rho),Z) \conv[\rho]{\infty} 0 .
\end{equation}
In particular, for {each nonnegative integer $r$} we have
\[\PPP{N_{W_\rho}(v_\rho) = {r}} \conv[\rho]{\infty}e^{-\tau}\cdot \frac{\tau^{r}}{{r}!}.\]
\end{theorem}

Now, for each $k\geq 1$, let $M^{(k)}_{W_\rho}$  be the $k$-th largest value of the inradius over all cells with incenter in $W_\rho$, and set $M^{(k)}_{W_\rho}=0$ if the number of such cells is smaller than $k$. The random variables $M^{(k)}_{W_\rho}$ are referred to as the \textit{order statistics}. In particular, when $k=1$, the random variable $M^{(1)}_{W_\rho}$ is the maximum of the inradii over all cells with incenter in $W_\rho$. As a corollary of Theorem \ref{Th:extremes}, we obtain the following result.
\begin{corollary}
\label{cor:extremes}
With the same assumptions as in Theorem \ref{Th:extremes} we have  for all $k\geq 1$
\[\PPP{M^{(k)}_{W_\rho}\leq v_\rho} \conv[\rho]{\infty} \sum_{r=0}^{k-1}e^{-\tau}\cdot \frac{\tau^r}{r!}.\]
\end{corollary}
This result is a direct consequence of Theorem \ref{Th:extremes} and of the fact that 
$M_{W_\rho}^{(k)}\leq v_\rho$ if and only if  $N_{W_\rho}(v_\rho)\leq k-1$. According to Corollary \ref{cor:extremes}, the maximum of inradii  belongs to the domain of attraction of a Gumbel distribution. Indeed, by taking $\tau=e^{-u}$ and $k=1$, we obtain 
\[\PPP{M_{W_\rho}^{(1)} \leq \frac{1}{2t}\log\rho + \frac{1}{2t}u}\conv[\rho]{\infty}e^{-e^{-u}},\] with $u\in \RR$.
This limit result  is classical in Extreme Value Theory according to  Gnedenko's theorem (see e.g. Theorem 1.1.3 in \cite{HF}). 

It is noticeable that the largest order statistics of the inradius in a STIT tessellation have asymptotically the same  distribution as the order statistics of the inradius in a Poisson line tessellation with the same intensity (see Theorem 1.1 (ii) in \cite{ChenHem}). This fact is not surprising in the sense that the typical cell of a STIT tessellation has the same distribution as the typical cell of a stationary Poisson line tessellation. However, Theorem \ref{Th:extremes} and Corollary \ref{cor:extremes} of the present paper are not consequences of Theorem 1.1 (ii) of \cite{ChenHem}. Indeed, dealing with extremes requires a specific treatment of correlations between cells and a mixing condition: the fact that two random tessellations have the same typical cell is not sufficient to ensure that the extremes have the same behaviour. Although Theorem  \ref{Th:extremes} and Theorem 1.1 (ii) of \cite{ChenHem} are both based on Poisson approximation, the methods which are used are different. In \cite{ChenHem}, the Poisson approximation is derived from the method of moments whereas, in our paper, this is based on the Chen-Stein method. We think that the method of moments is not suitable for STIT tessellations because they are based on a time process. A method similar to the one used in \cite{ChenHem} should lead to very technical computations in the context of STIT tessellation. \label{page:moment}

Our paper complements \cite{CC,Chen} where Poisson-Voronoi tessellations are treated. However, although  Theorem 1 in \cite{Chen} is a general result on extremes for random tessellations, the conditions of this theorem are too restrictive to be applied to STIT tessellations. Indeed, Theorem 1 in \cite{Chen} needs a so-called (FRC) condition which is not satisfied for a STIT tessellation. 

Our paper is organized as follows. In Section \ref{sec:preliminaries},  we set up the notation,
formally introduce the STIT tessellation processes  and recall several known results which will be used in the proof of Theorem \ref{Th:extremes}. In Section \ref{sec:technicalresults}, we state technical results which will be used to prove our main theorem. The proofs of  these technical lemmas and of Theorem \ref{Th:extremes}   are given in  Sections \ref{sec:prooftechnical} and \ref{sec:prooftheorem},  respectively. In Section \ref{sec:conclusion} we mention some potential extension of our main theorem. 

\section{Preliminaries}
\label{sec:preliminaries}

\subsection{Notation}
By $S^1$ we denote the unit sphere centered at the origin. For $c\in \RR^2$ and $R>0$, the set $B(c,R)$ denotes the topologically closed (Euclidean) disk with center $c\in \RR^2$ and radius $R>0$. For any subset $B\subset\RR^2$, the set  $B^\circ$ denotes the interior of $B$. If $B$ is a Borel subset, we recall that its area is denoted by $a(B)$. Further, $\oplus$ and $\ominus$ denote the Minkowski addition and subtraction of sets, respectively.

We denote by ${\cal H}$ the set of all lines in the plane. For any $B\subset \RR^2$,  we write  
\[[B]:= \{ H\in {\cal H}:\, H\cap B \not= \emptyset \}.\]
By $\Lambda$, we denote the measure on the set of lines ${\cal H}$ (endowed with the usual $\sigma$-algebra, see e.g. \cite{SW}), which is invariant under translations and rotations of the plane, with normalization $\Lambda([B(0,1)]) = 2$. In particular, for any convex polygon $z\subset \RR^2$, we have $\Lambda ([z])= perimeter (z)/\pi$. We shortly write $\mathrm{d}H$ instead of $\Lambda (\mathrm{d}H)$ in integrals over the set of lines. Thus, for all nonnegative measurable functions $f:\mathcal{H} \to \RR$, we have:
\begin{equation}
\label{eq:integralline} \int_{\mathcal{H}} f(H) \mathrm{d}H =\int_{(0,\infty )} \int_{S^1}  f(H(r,u)) \mathrm{d}r\, \sigma(\mathrm{d}u), 
\end{equation} 
where $H(r,u)$ is the line with distance $r>0$ from the origin and normal direction $u\in S^1$, and where $\sigma$ is the rotation invariant  measure on $S^1$ with normalization $\sigma(S^1)=2$. 

 Given three lines $H_1$, $H_2$ and $H_3$ in general position,  we use the notation $H_{1:3}:= (H_1,H_2,H_3)\in \mathcal{H}^3$, and we denote by $\Delta(H_{1:3})$  the unique triangle that can be formed by the intersection of  halfplanes induced by these lines, and  $B(H_{1:3})$, $c(H_{1:3})$, $R(H_{1:3})$ are the incircle, the incenter and the inradius of $\Delta(H_{1:3})$, respectively. In particular, $B(H_{1:3})=B(c(H_{1:3}), R(H_{1:3}))$. Similarly, if $s_1$, $s_2$, $s_3$ are three linear segments in $\RR^2$ in general position, we write $s_{1:3}:= (s_1,s_2,s_3)$. Correspondingly,
$B(s_{1:3})$, $c(s_{1:3})$, $R(s_{1:3})$ denote the incircle, the incenter and the inradius of the triangle formed by the three lines containing the three segments. {With a slight abuse of notation, we also write $\{s_{1:3}\}=\{s_1,s_2,s_3\}$.}

Let $Y_t$ be a STIT tessellation at time $t$. For any Borel subset $B\subset\RR^2$ and for any threshold $v\geq 0$, we write 
\begin{equation}\label{eq:notationmaxstit} M_B=\max_{z\in Y_t: c(z)\in B}R(z) \quad \text{and} \quad N_B(v)=\sum_{z\in Y_t:c(z)\in B}\ind{R(z)>v}.\end{equation} This notation $N_B(v)$ is consistent with \eqref{eq:defNWvrho}. 

Moreover, for any finite set $A$, we denote by $|A|$ the number of elements of $A$. For any real number $x$, we denote by $\lfloor x \rfloor$  the integer part of $x$. Given a probability space $(\Omega, \mathcal{E}, \PP)$, and two random variables $X$ and $Y$ defined on $\Omega$, the notation $X\stackrel{d}{=}Y$  indicates that $X$ and $Y$  have the same distribution. We  write $\EEE{Y|X}$  for the conditional expectation of $Y$ with respect to $X$.

\subsection{The STIT tessellation}
\label{STIT tessellation}
The STIT tessellation process in the $d$-dimensional Euclidean space was first introduced in \cite{NagelWeiss05}. Since then, in a series of publications several equivalent descriptions were given. In the present paper we consider planar stationary and isotropic STIT tessellations. Below we give  a short reminder for this particular case only.

We start with the construction of a tessellation process $(Y_{t,W},\, t\geq 0)$ in a bounded convex polygon $W$, referred to as a window. All the random variables that we consider are defined on some probability space $(\Omega, \mathcal{E},\PP)$. Let $\tau_0,\, \tau_1,\, \tau_2,\, \ldots $ be independent and identically distributed (i.i.d.) random variables, all exponentially distributed with parameter $1$.

\begin{enumerate}
\item[(i)] The initial state of the process is $Y_{0,W}=\{z_1\}:=\{ W\}$, and the random holding time in this state is $\tau_0/\Lambda ([W])$, i.e.  it is  exponentially distributed with parameter $\Lambda ([W])$. 

\item[(ii)]  At the end of the holding time, the window $W$ is divided by a random line $H_1$ with law  $(\Lambda ([W]))^{-1}\Lambda (\cdot \cap [W])$. This law is the probability distribution on $[W]$ generated by the restricted and normalized measure 
$\Lambda$. The new state of the STIT process is now $\{z_1,\, z_2\}$, where $z_1:= W\cap H_1^+$ and $z_2:= W\cap H_1^-$, and $H_1^+$ and $H_1^-$ are the two closed half-planes generated by $H_1$. Here it does not play a role which one is the positive or the negative half-plane. The random life times of  $z_1$ and $z_2$ are $\tau_1/\Lambda ([z_1])$ and  $\tau_2/\Lambda ([z_2])$, respectively.

\item[(iii)] Now, inductively, for $t>0$, assume that $Y_{t,W} =\{ z_{i_1},\ldots , z_{i_n}\}$. 
The life times of the cells are $\tau_{i_1}/\Lambda ([z_{i_1}]), \ldots , $ $\tau_{i_n}/\Lambda ([z_{i_n}])$, respectively. At the end of the life time of a cell $z_{i_j}$, this cell is divided by a random line $H_{i_j}$ with the law  $(\Lambda ([z_{i_j}]))^{-1}\Lambda (\cdot \cap [z_{i_j}])$, which is a probability distribution on $[z_{i_j}]$. Given the state of the tessellation process at the time of division, this line is conditionally  independent from all the other dividing lines used so far. The divided cell $z_{i_j}$ is deleted from the tessellation and is replaced by the two ``daughter'' cells $z_{i_j} \cap H_{i_j}^+$ and $z_{i_j} \cap H_{i_j}^-$. These cells are endowed with new indexes from $\NN$ which are not used before in this process.
\end{enumerate}

An essential property of the construction is that the distribution of the tessellation generated in a window $W$ is spatially consistent in the following sense. If $W$ and $W'$ are two convex polygons with $W\subset W'$ and $Y_{t,W}$, $Y_{t,W'}$  the respective random tessellations, then $Y_{t,W}\stackrel{d}{=} Y_{t,W'}\wedge W$  are identically distributed, where 
\[Y_{t,W'}\wedge W:=\{ z\cap W:\, z\in Y_{t,W},\, z\cap W^\circ\not= \emptyset \}\] 
is  the restriction of $Y_{t,W'}$ to $W$. This property yields the existence of a stationary random tessellation $Y_t$ of $\RR^2$ such that its restriction $Y_t\wedge W$ to any window $W$ has the same distribution as the constructed tessellation $Y_{t,W}$. Since the measure $\Lambda$ is invariant under rotation with the normalization $\Lambda([B(0,1)]) = 2$, the STIT tessellation $Y_t$ is isotropic and its intensity is $\gamma_t=\frac{t^2}{\pi}$ (see \cite{NagelWeiss08}). A ``global construction'', that does not refer to a bounded window, of the process $(Y_t,\, t>0)$ was given in \cite{mnw08a}, but this is much more involved.

Throughout this paper, we work with the stationary and isotropic STIT tessellation $Y_t$ for some $t>0$. Notice that a scaling property holds in the sense that
\begin{equation}\label{eq:scaling}
t\cdot Y_t := \{ t\cdot z:\, z \in Y_t\} \stackrel{d}{=}Y_1.
\end{equation}

The important constituents of planar STIT tessellations are the {\em maximal segments} (also referred to as I-segments), which are those segments generated by the intersection of a cell with its dividing line. In their relative interior, these maximal segments can contain endpoints of other maximal segments which appear later in the process. By $m_t$ we denote the set of all maximal segments of $Y_t$.  The union of the boundaries of the cells,  the so-called \textit{skeleton}, is denoted by $\partial Y_t$. It coincides with the union of all maximal segments.

\subsection{Poisson approximation}
It is clear that if $X(1),\ldots, X(\lfloor\rho\rfloor)$ are independent and identically distributed random variables with the same distribution as the inradius of the typical cell, i.e. with exponential distribution with parameter $2t$, then  $\sum_{i=1}^{ \lfloor \rho \rfloor}\ind{X(i)>v_\rho}$ converges in total variation to a Poisson random variable with parameter $\tau$, where $v_\rho$ is as in \eqref{eq:defvrho}. Theorem \ref{Th:extremes} establishes the same type of result, excepted that the random variables which we consider consist of a family of inradii of cells which are \textit{not independent}. The main difficulty in our work comes from the dependence between the cells, and consists in showing that the number of exceedances $N_{W_\rho}(v_\rho)$ has the same behaviour \textit{as if} we consider $\lfloor\rho\rfloor$ independent random variables, with the same distribution as the inradius of the typical cell.

In the same spirit as in \cite{Chen}, the main idea to derive a Poisson approximation of $N_{W_\rho}(v_\rho)$  is to apply a result due to Arratia, Goldstein and Gordon  \cite{AGG89}. Their result is based on the Chen-Stein method and gives an upper bound for the total variation distance between the distribution of a sum of Bernoulli random variables, and a Poisson distribution.

Let us recall the framework of their result. For  an arbitrary index set  $I$, and for $i\in I$, let $X_i$ be a Bernoulli random variable with $p_i=\PPP{X_i=1}=1-\PPP{X_i=0}$. For each $i,j\in I$, we  write $p_{ij}=\EEE{X_iX_j}$. Further, we let 
\[X:=\sum_{i\in I}X_i \quad \text{and} \quad \lambda:=\EEE{X}=\sum_{i\in I}p_i  , \quad \text{and assume that} \quad 0<\lambda <\infty .\]  
For each $i\in I$, fix a ``neighborhood'' $B_i\subset I$ with $i\in B_i$, and define
\[b_1:=\sum_{i\in I}\sum_{j\in B_i}p_ip_j, \quad b_2:=\sum_{i\in I}\sum_{i\neq j\in B_i}p_{ij}, \quad b_3 := \sum_{i\in I} \EEE{\left| \EEE{X_i-p_i\left| \sum_{j\in I\setminus B_i}X_j\right.}  \right|}. \]
Roughly, $b_1$ measures the neighborhood size, $b_2$ measures the expected number of neighbors of a given occurrence and $b_3$ measures the dependence between an event and the number of occurrences outside its neighborhood. We are now prepared to state a result on Poisson approximation (see Theorem 1 of \cite{AGG89}). 

\begin{proposition}(Arratia, Goldstein, Gordon)
\label{prop:AAG}
Let $Z$ be a Poisson random variable  with mean $\lambda \in (0,\infty)$. With the above notation and the  assumptions, we have 
\[d_{TV}(X,Z) \leq 2\, \left( (b_1+b_2)\cdot\frac{1-e^{-\lambda}}{\lambda} + b_3\cdot \min\{1,\, 1.4\,\lambda^{-1/2}\}  \right).\]
\end{proposition}

\section{Technical results}
\label{sec:technicalresults}

First notice, that for arbitrary $t>0$ the scaling property (\ref{eq:scaling}) of the STIT tessellation process yields that
$$
 \sum_{z\in Y_t:c(z)\in W_\rho}\ind{R(z)> v_\rho} \stackrel{d}{=} 
\sum_{z\in Y_1:c(z)\in t\, W_\rho}\ind{R(z)> t\, v_\rho}.
$$
Hence, it is sufficient to prove Theorem \ref{Th:extremes} for $t=1$, which will simplify the formulas. From now on, we will always set $t=1$, with the only exception in Lemma \ref{le:encapsulation}.

We adapt several arguments contained in \cite{Chen} to our context. The  difficulty compared to \cite{Chen} comes from the fact that a STIT tessellation only has a $\beta$-mixing  property, see \cite{MN}, whereas the tessellations considered in \cite{Chen} satisfy a finite range condition.  The main idea is to subdivide the window $W_\rho$ into small squares to be in the framework of Proposition \ref{prop:AAG}. We first introduce this subdivision. Then we establish several technical lemmas which deal with an asymptotic property and a local property, respectively. These technical lemmas will be applied to derive  Theorem \ref{Th:extremes} in the next section.

\subsection{A subdivision of the window}
\label{subsec:discretization}
Recall that $W_\rho=\sqrt{\pi \ \rho}\cdot [-\frac{1}{2},\frac{1}{2}]^2$ if $t=1$. When $\rho>e$, we subdivide the window $W_\rho$ into a set $V$ of sub-squares of equal size, where the number of these sub-squares is
 \begin{equation}\label{eq:subdivision} 
|V|=\left( \left\lfloor \sqrt{\frac{\pi \ \rho}{\log\log\rho}} \right\rfloor \right) ^2.
 \end{equation}
The sub-squares are indexed by  $\mathbf{i}:=(i_{1}, i_{2})\in \left[ 1,{\sqrt{|V|}}\right] ^{2}$, analogously to the order of indexing the elements of a matrix. With a slight abuse of notation, we identify a square with its index. 

In our proof we will often use the fact that there exists a $\rho_0(\tau)$ such that, for all $\rho>\rho_0(\tau)$, we have 
\begin{equation}\label{eq:rhotau}
v_\rho > \sqrt{2}\frac{\sqrt{\pi\ \rho}}{ \left\lfloor \sqrt{\frac{\pi \ \rho}{\log\log\rho}} \right\rfloor },
\end{equation}
where $v_\rho=v_\rho(\tau)$ is as in \eqref{eq:defvrho}. The right-hand side in the above equation is the length of the diagonal of a sub-square $\mathbf{i}$. Thus, for $\rho>\rho_0(\tau)$, there can be at most one incircle with center in $\mathbf{i}$ and radius larger than $v_\rho$.

The distance between sub-squares $\mathbf{i}$ and $\mathbf{j}$ is defined as $d(\mathbf{i},\mathbf{j}):=\max_{1\leq s\leq 2}|i_{s}-j_{s}|$. For $\mathbf{i}\in V$ and $r>0$ define the $r$-neighborhood of $\mathbf{i}$ as 
\begin{equation}\label{eq:defneighborhood}
S(\mathbf{i}, r) := \{\mathbf{j}\in V: d(\mathbf{i}, \mathbf{j})\leq r\}.
\end{equation}

The main idea to prove Theorem \ref{Th:extremes} is to apply Proposition \ref{prop:AAG} with $X_\mathbf{i}:= \ind{M_\mathbf{i}>v_\rho}$. We recall that $M_\mathbf{i}$ is the maximum of inradii over all cells with nucleus in $\mathbf{i}$, i.e. $M_\mathbf{i}=\max_{z\in {Y_1}: c(z)\in \mathbf{i}}R(z)$. The sets $B_i$ of Proposition \ref{prop:AAG} are replaced by the neighborhoods $S(\mathbf{i},\rho^{\beta /2} )$ for some $\beta \in (0,1)$. Furthermore, for all sub-squares $\mathbf{i}, \mathbf{j}\in V$, we denote
\[p_\mathbf{i} := \PPP{M_\mathbf{i}>v_\rho} \quad \text{and} \quad p_{\mathbf{ij}} := \PPP{M_\mathbf{i}>v_\rho, M_\mathbf{j}>v_\rho}.\]
We also let
\begin{equation}
\label{eq:defb1b2b3}
 b_1:=\sum_{\mathbf{i}\in V}\sum_{\mathbf{j}\in S(\mathbf{i},\rho^{\beta /2} )}p_\mathbf{i}p_\mathbf{j}, \quad 
b_2:=\sum_{\mathbf{i}\in V}\sum_{\mathbf{i}\neq \mathbf{j}\in S(\mathbf{i},\rho^{\beta /2} )}p_{\mathbf{ij}}, \quad 
b_3 := \sum_{\mathbf{i}\in V} {\EEE{\left| \PPP{ M_\mathbf{i}>v_\rho\left| \sum_{\mathbf{j}\not\in S(\mathbf{i},\rho^{\beta /2} )} \ind{M_\mathbf{j}>v_\rho}\right.}- p_\mathbf{i}  \right|}  }.
\end{equation}
In Section \ref{sec:prooftheorem}, we will show that $b_1$, $b_2$ and $b_3$ converge to 0 as $\rho$ goes to infinity, and thus an application of Proposition \ref{prop:AAG} yields the proof of Theorem \ref{Th:extremes}.

\subsection{Technical lemmas} \label{techlemmas}
We start with some  technical lemmas which will be used to prove that $b_2$ and $b_3$ converge to 0 as $\rho$ goes to infinity.

\paragraph{Lemmas concerning an upper bound for $b_2$}
The arguments showing that $b_2$ converges to 0 are mainly inspired by \cite{ChenHem}, and they rely on a ``local condition'', i.e. on an upper bound for the probability that the incircles of two cells with centers in a short distance simultaneously exceed the threshold $v_\rho$. The following  lemma  provides an upper bound of the probability that $\partial Y_1$ does not intersect a union of disks, and it is an adaptation of Lemma 4.4 of \cite{ChenHem} in the context of STIT tessellations. 
\begin{lemma}
\label{Le:twoballs}
There exists a constant $\eta >0$ such that for all pairs of disks
 $B_1=B(c_1,r)$ and $B_2=B(c_2,r)$  with the same radius $r>0$ and with centers $c_1,c_2\in \RR^2$ satisfying $||c_1-c_2||\geq 2r$ we have
\[\PPP{\partial Y_1\cap (B_1\cup B_2)=\emptyset} \leq \eta \cdot e^{-2(1+\frac{2}{\pi}) r}.\]
\end{lemma}
The proof is given in Subsection \ref{sec:prooftechnical}.
This lemma is one of the key ingredients to prove that $b_2$ converges to 0. 
Indeed, together with some non-trivial computations which are performed in Section \ref{sec:prooftheorem},  Lemma \ref{Le:twoballs}  yields that, with high probability, the inradii of two cells which are close enough (in the sense that the incenter of one of them belongs to some sub-square $\mathbf{i}$, and the other belongs to some sub-square $\mathbf{j}\in S(\mathbf{i}, \rho^{\beta})$) cannot simultaneously exceed the threshold $v_\rho$. 

The following lemmas  are adaptations of Lemma A.1 (i) and Lemma A.1 (ii) in \cite{ChenHem} respectively.
\begin{lemma}\label{lem:G(H_1)}
There exists a constant $\eta >0$ such that for all lines $H_1\in \mathcal{H}$ and all $R >0$  we have
 \begin{equation*}
        G(H_1):=  \int_{\mathcal{H}^2}\ind{c(H_{1:3})\in B(0,R)}\mathrm{d}H_{2:3} \leq     \eta \,  {R} ^2 .
    \end{equation*}
\end{lemma} 
The following result is similar to Lemma \ref{lem:G(H_1)}, but this time two lines are fixed.
\begin{lemma}\label{lem:G(H_2)}
There exists a constant $\eta >0$ such that for all pairs of lines $H_1, H_2 \in \mathcal{H}$ and all $R >0$  we have \begin{equation*}
        G(H_1, H_2):=  \int_{\mathcal{H}}\ind{c(H_{1:3})\in B(0,R)}\mathrm{d}H_{3}\leq    \eta \,  R .
    \end{equation*}
\end{lemma}

\paragraph{Lemmas concerning an upper bound for $b_3$}
The arguments showing that $b_3$ converges to 0 are mainly inspired by \cite{MN}. The key argument is a mixing property of the STIT tessellations. However, the general upper bound for the $\beta$-mixing coefficient provided in \cite{MN} is not sufficient for our purposes. Therefore, a more specific treatment of rare events is developed. 

To do it, let $\mathcal{C}$ denote the set of all compact convex subsets of $\RR^2$. Furthermore, let $\sigma (\mathcal{I})$ be the  $\sigma$-algebra generated by the family  $\mathcal{I}$ of sets. Denote by $\mathcal{T}$ the set of all tessellations of $\RR^2$. We endow $\mathcal{T}$ with the Borel $\sigma$-algebra $\mathcal{B}(\mathcal{T})$ of the Fell topology, namely
\[\mathcal{B}(\mathcal{T}):=\sigma \left( \{ \{ T\in\mathcal{T} :\, \partial T\cap C=\emptyset \}     :\, C\in\mathcal{C}  \}   \right).\] \label{page:tessellation}
Moreover, for some compact convex set $K\subset \RR^2$, we define
\begin{align*}
& \mathcal{B}(\mathcal{T}_{K}):=\sigma \left( \{ \{ T\in\mathcal{T} :\, \partial T\cap C=\emptyset \}     :\, C\subset K,\, C\in\mathcal{C}  \}   \right),\\
& \mathcal{B}(\mathcal{T}_{K^c}):=\sigma \left( \{ \{ T\in\mathcal{T} :\, \partial T\cap C=\emptyset \}     :\, C\subset K^c,\, C\in\mathcal{C}  \}   \right).
\end{align*}

Let $K',K$ be two  convex polygons such that $0\in K'\subset K$. In \cite{MN2} the concept of \textit{encapsulation} was  introduced. It means that there is a state of the STIT process $\mathbf{Y}=(Y_t,t>0)$  such that all facets of $K'$ are separated from the facets of $K$ by facets of the tessellation before the interior of $K'$ is divided by a facet of the tessellation. Formally, denoting the $0$-cell by $C_t^0$, i.e.  the cell of $Y_t$ that contains the origin, we define the encapsulation time as
\[S(K,K') := \inf\{t>0: K'\subset C_t^0\subset K^\circ\},\] 
with the convention $\inf\emptyset=\infty$. 
The following lemma is an adaptation of Lemma 6.4 in \cite{MN}. 

\begin{lemma}
\label{le:encapsulation}
Let $\mathbf{Y}=(Y_t, t>0)$ be a STIT process.  Let $0\in K'\subset K$ be two convex polygons,  $\mathcal{E}\in \mathcal{B}(\mathcal{T}_{K^c})$ and $\mathcal{A}\in \mathcal{B}(\mathcal{T}_{K'})$. Then for all $0<s<t$ we have
\begin{align*}
& \left|\PPP{Y_t\in \mathcal{E}|Y_t\in \mathcal{A}} - \PPP{Y_t\in \mathcal{E}}\right| \\
\leq  & \PPP{Y_t\in \mathcal{E}|Y_t\in \mathcal{A}} -  \frac{\PPP{Y_{t-s}\in \mathcal{A}}}{\PPP{Y_t\in \mathcal{A}}}\cdot \PPP{Y_t\in \mathcal{E}, S(K,K')<s, Y_s\wedge K'=K'}\\
& + \left| \frac{\PPP{Y_{t-s}\in \mathcal{A}}}{\PPP{Y_t\in \mathcal{A}}} - 1  \right|\PPP{Y_t\in \mathcal{E}} + \frac{\PPP{Y_{t-s}\in \mathcal{A}}}{\PPP{Y_t\in \mathcal{A}}}\cdot \PPP{\{Y_t\in \mathcal{E}\}\cap \{S(K,K')<s, Y_s\wedge K'=K'\}^c}.
\end{align*}
\end{lemma}

Now we construct particular squares $K'\subset K$ which are tailored for our purposes.
Recall that $V$ is a subdivision of the window $W_\rho$, $\rho>0$, into small squares, and that $S(\mathbf{i},r)$ denotes the $r$-neighborhood of a square $\mathbf{i}\in V$ (see Section \ref{subsec:discretization}). By $\mathbf{S}_0$ we denote the square centered at the origin with side length $2 v_\rho^2$, and 
\begin{equation}\label{eq:defCi} 
C(\mathbf{i}):= W_\rho \setminus \bigcup_{\mathbf{j}\in S(\mathbf{i},\rho^{\beta /2} )} \mathbf{j}
\end{equation} 
for the complement with respect to $W_\rho$ of the neighborhood of  $\mathbf{i}\in V$.

Now, because in stationary and isotropic STIT tessellations there are a.s. no cells with pairs of parallel sides, we can restrict the considerations to the measurable set  $\mathcal{T}'\subset\mathcal{T}$  of all tessellations which do not contain cells with pairs of parallel sides. This simplifies  the study of incircles of the cells.
 For any $T\in \mathcal{T}'$, for any Borel subset $B\subset\RR^2$, and for any threshold $v\geq 0$, let 
\begin{equation*}
M_{B}^{(T)}:=\max_{z\in T: c(z)\in B}R(z) \quad \text{and} \quad  N_{B}^{(T)}(v):=\sum_{z\in T: c(z)\in B}\ind{R(z)>v}. \end{equation*} 
 Notice that  $M_{(C(\mathbf{i}))^\circ}^{(T)}$ and $N_{(C(\mathbf{i}))^\circ}^{(T)}(v_\rho)$ are the maximum of inradii and the number of exceedances \textit{outside} a neighborhood of $\mathbf{i}$. The following lemma provides the background to apply Lemma \ref{le:encapsulation}.

\begin{lemma}
\label{le:measurable}
Let $v_\rho$ be as in \eqref{eq:defvrho}, with $t=1$ and $\rho>\rho_0(\tau) $ (see (\ref{eq:rhotau})), and $\mathbf{i}\in V$ be fixed. Let $K=\left(\bigcup_{\mathbf{j}\in S(\mathbf{i},\rho^{\beta /2} )} \mathbf{j}\right) \ominus \mathbf{S}_0$ and $K'=\mathbf{i}\oplus \mathbf{S}_0$.
Then 
\begin{enumerate}[(i)]
\item $\left\{T\in \mathcal{T}': v_\rho <M_{\mathbf{i}^\circ}^{(T)}\leq v_\rho^2 \right\} \in \mathcal{B}({\cal T}_{K'}).$
\item $\left\{T\in \mathcal{T}': M_{(C(\mathbf{i}))^\circ}^{(T)}\leq v_\rho^2 \quad \text {and}\quad N_{(C(\mathbf{i}))^\circ}^{(T)}(v_\rho) = k\right\}\in \mathcal{B}({\cal T}_{K^c})$ for all $k\in \NN_0$.
\end{enumerate}
\end{lemma}

\section{Proofs of technical lemmas}
\label{sec:prooftechnical}
\begin{prooft}{Lemma \ref{Le:twoballs}}
Let $B_1=B(c_1,r)$ and $B_2=B(c_2,r)$ be two disks  with the same radius $r>0$ and with centers $c_1,c_2\in\RR^2$. Denote by  $[B_1|B_2]$ the set of  lines which separate $B_1$ and $B_2$, and write  $\text{conv}(B_1\cup B_2)$ for the convex hull of $B_1\cup B_2$. According to Corollary 1 in \cite{NagelWeiss05}, the probability $\PPP{\partial Y_1\cap (B_1\cup B_2) = \emptyset}$ can be written as follows:
\begin{itemize}
\item if $\Lambda([B_1])+\Lambda([B_2])\neq \Lambda([\text{conv}(B_1\cup B_2)])$, we have 
\begin{multline}\label{eq:convexhull} 
\PPP{\partial Y_1\cap (B_1\cup B_2) = \emptyset} = e^{-\Lambda([\text{conv} (B_1\cup B_2)])}\\ + \Lambda([B_1|B_2])\cdot \frac{e^{-\Lambda([\text{conv}(B_1\cup B_2)])} - e^{-\Lambda([B_1])}e^{-\Lambda([B_2])}}{\Lambda([B_1])+\Lambda([B_2])-\Lambda([\text{conv}(B_1\cup B_2)])},
\end{multline}
\item otherwise, we have 
\begin{equation} \label{eq:convexhull2}
\PPP{\partial Y_1\cap (B_1\cup B_2) = \emptyset} = e^{-\Lambda([\text{conv}(B_1\cup B_2)])}+ \Lambda([B_1|B_2])\, e^{-\Lambda([B_1])}e^{-\Lambda([B_2])}.
\end{equation}
\end{itemize}

Now, let $d:=||c_2-c_1||\geq 2r$. It is clear that
\[\Lambda([B_1])=\Lambda([B_2]) = 2 r \quad \text{and} \quad \Lambda([\text{conv}(B_1\cup B_2])) = 2 r+2d/\pi .\]
Moreover, because $[B_1|B_2] \subseteq  [\overline{c_1 c_2}]$, we obtain $\Lambda([B_1|B_2]) \leq  \Lambda([\overline{c_1 c_2}])= \frac{2d}{\pi}$,  where $\overline{c_1 c_2}$ denotes the linear segment connecting $c_1$ and $c_2$. 

We discuss the two cases \eqref{eq:convexhull} and \eqref{eq:convexhull2}. First, assume that $d\neq \pi r$ and hence $\Lambda([B_1])+\Lambda([B_2])\neq \Lambda([\text{conv}(B_1\cup B_2)])$. Thus, according to \eqref{eq:convexhull}, we have
\[\PPP{\partial Y_1\cap (B_1\cup B_2) = \emptyset} \leq e^{-(2 r +2d/\pi)} + \frac{2d}{\pi}\cdot \frac{e^{-(2 r +2d/\pi )} - e^{-4 r} }{2 r-2d/\pi}.\]
Now, we provide  an upper bound, which is independent of $d$,  for the right-hand side.
Taking $u:=\frac{2d}{\pi}-2 r$, we have $u\geq (\frac{4}{\pi}-2)r$ and
\begin{align*}
\PPP{\partial Y_1\cap (B_1\cup B_2) = \emptyset} &   \leq e^{-(u+ 4 r)} + (u+2r)\, \frac{e^{-(u+4 r)} - e^{-4 r} }{-u} \\
&=  e^{-4r}\,  \left(e^{-u} +  (u+2r) \cdot \frac{e^{-u} - 1 }{-u}  \right)\\
&=e^{-4r}\, \left(1+2r\cdot \frac{1-e^{-u}}{u}\right).
\end{align*}
Notice that the function $f$ defined as $f(u)=\frac{1-e^{-u}}{u}$ if $u\neq 0$, and $f(0):=1$, is positive, continuous and strictly decreasing on $\RR$. 
Since  $u\geq (\frac{4}{\pi}-2)r$, we have $\frac{1-e^{-u}}{u}\leq \frac{1-e^{- (\frac{4}{\pi}-2)r}}{ (\frac{4}{\pi}-2)r}$.
Hence, 
\begin{align*}
\PPP{\partial Y_1\cap (B_1\cup B_2) = \emptyset} & \leq e^{-4r}\left(1 +\frac{2}{2-\frac{4}{\pi}}\left(e^{(2-\frac{4}{\pi})r}-1 \right) \right)\\
& \leq \frac{2}{2-\frac{4}{\pi}} \, e^{-2(1+\frac{2}{\pi}) r}.
\end{align*}

In the case $d=\pi r$, the equation \eqref{eq:convexhull2} yields 
\begin{align*}
\PPP{\partial Y_1\cap (B_1\cup B_2) = \emptyset} &
\leq e^{-(2 r +2d/\pi)} + \frac{2d}{\pi}\,  e^{-4 r} \\
& =  (1+2r)\, e^{-4r}\\
& = o(e^{-2(1+\frac{2}{\pi}) r}),
\end{align*}
which proves the assertion.
\end{prooft}

\bigskip

\begin{prooft}{Lemma \ref{lem:G(H_1)}}
 Let $H_1\in \mathcal{H}$ be fixed, with normal direction $u_1\in S^1$. According to \eqref{eq:integralline}, we have 
\[G(H_1)  =  \int_{(S^1)^2} \int_{\RR_+^2} \ind{c(H_1,H(r_2,u_2), H(r_3,u_3))\in B(0,R)}
 \mathrm{d}r_{2:3}\, \sigma^{\otimes 2}(\mathrm{d} u_{2:3}).\]

Now, let $(u_2,u_3)\in (S^1)^2$ and $(r_2,r_3)\in (0,\infty )^2$ such that $H_{1:3}=(H_1,H_2,H_3)$ are in general position, with $H_2=H(r_2,u_2)$ and $H_3=H(r_3,u_3)$. For short, we write $c:=c(H_{1:3})$. Let $<\cdot,\cdot >$ be the scalar product and identify the points $c$ and $u_i$ with their corresponding vectors. Furthermore, let $r(c,H)$ denote the distance of a line $H$  to a point $c\in\RR^2$. Elementary geometry yields $r(c,H_i)=r(0,H_i) - <u_i,c>$ , $i=1,2,3$. Hence, for $i=2,3$, we have \begin{equation}\label{eq:scalarproduct}r_i=r(0,H_i)=r(0,H_1) - <u_1,c> + <u_i,c>.\end{equation} Now,  for each $u_2,u_3\in S^1$, consider the change of variables $c\mapsto (r_2,r_3)$, where $r_2,r_3$ are two positive numbers such that $c(H_1,H(r_2,u_2), H(r_3,u_3))=c$. According to  \eqref{eq:scalarproduct}, this map is linear and it is one-to-one. Its Jacobian does not depend on $c$, but only on $u_i$, $i=1,2,3$, and is be bounded by a constant. This shows that $G(H_1)$ is lower than a constant multiplied  by $\int_{(0,\infty )^2}\ind{c\in B(0,R)}\mathrm{d}c=\pi R^2$, which concludes the proof of Lemma \ref{lem:G(H_1)}.
\end{prooft}

\bigskip

\begin{prooft}{ Lemma \ref{lem:G(H_2)}}
This will be sketched only because it relies on the same idea as the proof of Lemma \ref{lem:G(H_1)}, see also the proof of Lemma A.1, (ii) of \cite{ChenHem}. Let $H_1, H_2 \in \mathcal{H}$ be fixed. For each line $H_3$, the incenter $c(H_{1:3})$ belongs to one of the two bisecting lines associated with $H_1$ and $H_2$. By considering a linear change of variables whose  Jacobian only depends on the normal vectors, we see that $G(H_1,H_2)$ is smaller than the diameter of $B(0,R)$ and thus $G(H_1,H_2) \leq \eta\, R$ for some constant $\eta$.
\end{prooft}

\bigskip

\begin{prooft}{Lemma \ref{le:encapsulation}}
First we write for $0<s<t$
\begin{multline*}
\left|\PPP{Y_t\in \mathcal{E}|Y_t\in \mathcal{A}} - \PPP{Y_t\in \mathcal{E}}\right|\leq \left|  \PPP{Y_t\in \mathcal{E}|Y_t\in \mathcal{A}} - \PPP{Y_t\in \mathcal{E}, S(K,K')<s, Y_s\wedge K'=K'| Y_t\in \mathcal{A}} \right|\\
+  \left|\PPP{Y_t\in \mathcal{E}, S(K,K')<s, Y_s\wedge K'=K'| Y_t\in \mathcal{A}} - \PPP{Y_t\in \mathcal{E}}\right|.
\end{multline*}
According to Lemma 6.2 of \cite{MN}, we can write the first term on the right-hand side as
\begin{multline*}
\PPP{Y_t\in \mathcal{E}|Y_t\in \mathcal{A}} - \PPP{Y_t\in \mathcal{E}, S(K,K')<s, Y_s\wedge K'=K'| Y_t\in \mathcal{A}}\\
 = \PPP{Y_t\in \mathcal{E}|Y_t\in \mathcal{A}} -  \PPP{Y_t\in \mathcal{E}, S(K,K')<s, Y_s\wedge K'=K'}\cdot \frac{\PPP{Y_{t-s}\in \mathcal{A}}}{\PPP{Y_t\in \mathcal{A}}}.
\end{multline*}

For the second term, we apply Lemma 6.2 of \cite{MN} again and write
\begin{align*}
 &\left|\PPP{Y_t\in \mathcal{E}, S(K,K')<s, Y_s\wedge K'=K'| Y_t\in \mathcal{A}} - \PPP{Y_t\in \mathcal{E}}\right| \\
\leq &\left|\PPP{Y_t\in \mathcal{E}, S(K,K')<s, Y_s\wedge K'=K'}\cdot \frac{\PPP{Y_{t-s}\in \mathcal{A}}}{\PPP{Y_t\in \mathcal{A}}} - \PPP{Y_t\in \mathcal{E}}\cdot \frac{\PPP{Y_{t-s}\in \mathcal{A}}}{\PPP{Y_t\in \mathcal{A}}} \right| \\
& + \left| \PPP{Y_t\in \mathcal{E}}\cdot \frac{\PPP{Y_{t-s}\in \mathcal{A}}}{\PPP{Y_t\in \mathcal{A}}} -    \PPP{Y_t\in \mathcal{E}} \right|
\\
=&  \frac{\PPP{Y_{t-s}\in \mathcal{A}}}{\PPP{Y_t\in \mathcal{A}}}\PPP{\{Y_t\in \mathcal{E}\}\cap \{S(K,K')<s, Y_s\wedge K'=K'\}^c}\\
 & + \left| \frac{\PPP{Y_{t-s}\in \mathcal{A}}}{\PPP{Y_t\in \mathcal{A}}}  - 1 \right|\PPP{Y_t\in \mathcal{E}},
\end{align*}
which concludes the proof of Lemma \ref{le:encapsulation}. 
\end{prooft}

\bigskip

\begin{prooft}{Lemma \ref{le:measurable}}
First, we prove (i). To do it, for each $x\in \RR^2$, $T\in \mathcal{T}'$, we denote by $d_T(x)$ the distance from the point $x$ to the skeleton $\partial T$ of the tessellation $T$. For a fixed $T$ the function $d_T$ is continuous, piecewise linear, and its local maxima are the incenters of the cells.

For all $x\in \mathbf{i}$, and all $r\leq v_\rho^2$, we have  $B(x,r)\subset K'$, and therefore
\[\{d(x)>r\} :=  \{ T\in {\mathcal{T}'}:\, d_T(x) >r \}  =\{ T\in {\mathcal{T}'}:\, B(x,r) \cap \partial T= \emptyset  \} \in \mathcal{B}(\mathcal{T}_{K'}).\]

  Let $\QQ$ denote the set of rational numbers. The idea for the formal proof given below is the following. Choose an $n_0\in \NN_0$,  large enough such that  a \textit{local} maximum $c(z)$ of $d_T$ inside a cell $z$ is focused on. Then this local maximum can be approximated by a sequence of points $x_n$, $n>n_0$, with rational coordinates. A point $x_n$ can be chosen in a distance $\frac{1}{2n}$ from $c(z)$, not arbitrarily, but in a direction from $c(z)$ where the slope of the function $d_T$ is smaller or equal than in the other directions.
Notice that with this method the incenters which are located in the interior $\mathbf{i}^\circ$ of $\mathbf{i}$ can be identified only; for incenters on the boundary of $\mathbf{i}$  we would need information outside $K'$. In a stationary random tessellation there is almost surely no incenter on the boundary of $\mathbf{i}$.

For $\rho>\rho_0(\tau)$ and $T\in \mathcal{T}'$, we have the following logical equivalences:
\begin{align*}
     & v_\rho <M_{\mathbf{i}^\circ}^{(T)}\leq v_\rho^2 \\
\iff & \exists z\in T:\, c(z)\in \mathbf{i}^\circ ,\, v_\rho <M_{\mathbf{i}^\circ}^{(T)}\leq v_\rho^2\\
\iff & \exists n_0\in \NN \, \forall n\geq n_0 \, \exists x_n \in  (\mathbf{i}^\circ\cap \QQ^2) \, \forall x\in 
          \left(B\left( x_n,\frac{1}{n_0}\right)\setminus B\left( x_n, \frac{1}{n}\right)\right)\cap \QQ^2 :\\
      & d_T(x)< d_T(x_n) \mbox{ and } v_\rho <d_T(x_n) \leq v_\rho^2.
\end{align*}

Thus,  $\{T\in \mathcal{T}': v_\rho <M_{\mathbf{i}^\circ}^{(T)}\leq v_\rho^2 \}$ can be represented as a union, intersection or complement of countably many events of the type $\{d(x)>r\}$, with $x\in \mathbf{i}^\circ$ and $r\leq v_\rho^2$, and hence it is in
$\mathcal{B}({\cal T}_{K'})$. This concludes the proof of (i).

The proof of (ii) relies on the same idea as the proof of (i). Observe that for all $x\in (C(\mathbf{i}))^\circ$, and all $r\leq v_\rho^2$, we have that $B(x,r)\subset K^c$, and therefore $\{d(x)>r\}   \in \mathcal{B}(\mathcal{T}_{K^c})$. Moreover, for any pair of points $x_1,\, x_2 \in  (C(\mathbf{i}))^\circ$, the linear segment $\overline{x_1 x_2}\subset K^c$ and therefore
$\{ T\in \mathcal{T}': \overline{x_1 x_2} \cap \partial T \not= \emptyset\} \in \mathcal{B}(\mathcal{T}_{K^c})$.
The condition $ \overline{x_{n,i} x_{n,j}} \cap \partial T \not= \emptyset $ ensures that $x_{n,i}$ and $x_{n,j}$ are located in different cells.
For all $T\in \mathcal{T}'$ we have
\begin{equation*}  M_{(C(\mathbf{i}))^\circ}^{(T)}\leq v_\rho^2 
\iff  \forall x\in (C(\mathbf{i}))^\circ \cap \QQ^2 :\, B(x,v_\rho^2)\cap \partial T \not= \emptyset .
\end{equation*}
Furthermore, for all $T\in \mathcal{T}'$,
\begin{align*}
     & N_{(C(\mathbf{i}))^\circ}^{(T)} \geq k \\
\iff & \exists n_0\in \NN \, \forall n\geq n_0 \, \exists (x_{n,1},\ldots, x_{n,k}) \in 	\left( (C(\mathbf{i}))^\circ \cap \QQ^2 \right)^k : \\
     & \left(\forall 1\leq i<j\leq k : \overline{x_{n,i} x_{n,j}} \cap \partial T \not= \emptyset \right) \mbox{ and } \\
		  & \forall j\in \{1,\ldots ,k\}  \forall x\in 
          \left(B\left( x_{n,j},\frac{1}{n_0}\right)\setminus B\left( x_{n,j}, \frac{1}{n}\right)\right)\cap \QQ^2 :
          d_T(x)< d_T(x_{n,j}) \mbox{ and } v_\rho <d_T(x_n) .
\end{align*}

 Thus, the set $\left\{T\in \mathcal{T}':  N_{(C(\mathbf{i}))^\circ}^{(T)} \geq k \quad \text{and} \quad M_{(C(\mathbf{i}))^\circ}^{(T)}\leq v_\rho^2\right\}$ can be represented as a union, intersection or complement of countably many events of the type $\{d(x)>r\}$, with $x\in \mathbf{i}$ and $r\leq v_\rho^2$, and hence it is in $\mathcal{B}({\cal T}_{K^c})$.
\end{prooft}

\section{Proof of Theorem \ref{Th:extremes}}
\label{sec:prooftheorem}

Recall that it is sufficient to prove the theorem for the case $t=1$ and that for $\rho>\rho_0(\tau)$ (see \eqref{eq:rhotau}) there can be at most one incircle with center in $\mathbf{i}\in V$ and radius larger than $v_\rho$, so that
\begin{equation}\label{eq:N=X}
N_{W_\rho}(v_\rho) = \sum_{\mathbf{i}\in V}\ind{M_\mathbf{i}>v_\rho}.
\end{equation}

Now, according to \eqref{Campbellmark} and \eqref{eq:N=X},  for all $\rho>\rho_0(\tau)$, we obtain that 
\begin{equation}\label{eq:sum pi}
\sum_{\mathbf{i}\in V}p_\mathbf{i} 
 = \EEE{\sum_{\mathbf{i}\in V}\ind{M_\mathbf{i}>v_\rho}}
 =  \EEE{N_{W_\rho}(v_\rho)}
 = \rho \ \PPP{R(\mathcal{Z})>v_\rho}
 =\tau.
\end{equation}
Thus, to prove Theorem \ref{Th:extremes}, it suffices to show that $\sum_{\mathbf{i}\in V}\ind{M_\mathbf{i}>v_\rho}$ converges in total variation to a Poisson random variable $Z$ with mean $\EEE{Z}=\tau$. According to Proposition \ref{prop:AAG}, it is sufficient to prove  that $b_1$, $b_2$ and $b_3$,  defined in  \eqref{eq:defb1b2b3}, converge to 0 as $\rho$ goes to infinity. These proofs are given in Sections   \ref{b1}, \ref{sec:localdependence} and \ref{subsection:mixing}, respectively.

\subsection{Proof of the assertion that $b_1$ converges to 0 for $\rho \to \infty$}\label{b1}
Because the STIT tessellation $Y_1$ is  stationary, the value of $p_\mathbf{i}$ does not depend on $\mathbf{i}$. Therefore, for $\mathbf{i}\in V$, we have
\[b_1\le |V|\cdot \sup_{\mathbf{i}\in V}|S(\mathbf{i}, \rho^{\beta /2})|\cdot p_\mathbf{i}^2
\leq  |V|\cdot \left( 2\lfloor \rho^{\beta/2}\rfloor+1  \right)^2\cdot p_\mathbf{i}^2 .\]
Now, for $\rho > \rho_0 (\tau)$, we observe that 
 \begin{equation}\label{eq:pi}
p_\mathbf{i} = \PPP{M_\mathbf{i}>v_\rho} = \EEE{\sum_{z\in Y_1: c(z)\in \mathbf{i}}\ind{R(z)>v_\rho}}= a(\mathbf{i})\ \gamma_1 \ \PPP{R(\mathcal{Z})>v_\rho} =  a(\mathbf{i})\ \gamma_1 e^{-2v_\rho},
\end{equation}  
where the last two equalities come from \eqref{eq:deftypical} and \eqref{eq:exponential}. Since the area of $\mathbf{i}$ is $a(\mathbf{i}) = \frac{ \pi\rho}{|V|}$, it follows from \eqref{eq:defvrho} that 
\[b_1\leq  \frac{\tau^2}{|V|} \, \left( 2\rho^{\beta/2}+1  \right)^2.\]

According to \eqref{eq:subdivision}, this implies that $b_1 = O\left(\rho^{-(1-\beta)}\log\log\rho\right)$.  Since $\beta<1$, this proves that $b_1$ converges to 0 as $\rho$ goes to infinity.

\subsection{Proof of the assertion that $b_2$ converges to 0 for $\rho \to \infty$}
\label{sec:localdependence}
Let $\mathbf{i}, \mathbf{j}\in V$ be fixed, with $\mathbf{i}\not= \mathbf{j}$. In the same spirit as in \eqref{eq:N=X}, for $\rho > \rho_0 (\tau)$, there are at most two incircles with radii larger than $v_\rho$ and centers in $\mathbf{i}$ and $\mathbf{j}$ respectively. Therefore
\begin{equation} 
\label{eq:defpij}p_{\mathbf{ij}} = \PPP{M_\mathbf{i}>v_\rho, M_\mathbf{j}>v_\rho} = \EEE{\sum_{(z_1,z_2)\in Y_1^2: (c(z_1), c(z_2))\in \mathbf{i}\times \mathbf{j}}\ind{R(z_1)>v_\rho}\ind{R(z_2)>v_\rho}}.
\end{equation}
To deal with the right-hand side, we consider three cases, regarding the number of tangential (to the incircle) maximal segments that are common to both cells.   By a common tangential maximal segment (CTS) of the cells $z_1$ and $z_2$ we mean a maximal segment $s\in m_1$ which is tangential to the incircles of $z_1$ and $z_2$ simultaneously, i.e. $s\cap B(z_1)\not= \emptyset$ and $s\cap B(z_2)\not= \emptyset$.  Notice that the cells $z_1$ and $z_2$ can have a common segment which is not a CTS.
Thus we write  $p_{\mathbf{ij}}$ as
\[p_{\mathbf{ij}}=p_{\mathbf{ij}}(0)+p_{\mathbf{ij}}(1)+p_{\mathbf{ij}}(2),\]
where for  $0\leq k\leq 2$,
\[p_{\mathbf{ij}}(k):=\EEE{\sum_{(z_1,z_2)\in Y_1^2: (c(z_1), c(z_2))\in \mathbf{i}\times \mathbf{j}}\ind{R(z_1)>v_\rho}\ind{R(z_2)>v_\rho}\ind{\text{($z_1,z_2$ have $k$ CTS)}}}\]
We prove below that $p_{\mathbf{ij}}(k)= O\left({(\log\log\rho)^2\cdot \log\rho}\cdot \rho^{-2(1+\epsilon)}\right)$ for each $k=0,1,2$ 
and   $0< \epsilon{\leq}\frac{2}{\pi}$. We begin with the case $k=0$. To do it, we will use two  formulas of Stochastic Geometry and Integral Geometry respectively, namely a Mecke-type formula for STIT tessellations and a Blaschke-Petkantschin type change of variables. The cases $k=1$ and $k=2$ will be dealt in a similar way, but also applying Lemmas \ref{lem:G(H_1)} and \ref{lem:G(H_2)}, respectively.

\subsubsection{Case $k=0$}		
We rewrite $p_{\mathbf{ij}}(0)$ in terms of maximal segments (instead of cells) of $Y_1$. This yields
\begin{equation*}
p_{\mathbf{ij}}(0)= \frac{1}{3! \cdot 3!}\EEE{\sum_{s_{1:3}\in m_1^3}\ \sum_{q_{1:3}\in m_1^3 \setminus \{ s_{1:3}\}}\ind{R(s_{1:3})>v_\rho, R(q_{1:3})>v_\rho} 
\ind{\partial Y_1\cap (B^o(s_{1:3})\cup B^o(q_{1:3}))=\emptyset}\ind{(c(s_{1:3}), c(q_{1:3}))\in \mathbf{i}\times \mathbf{j}}}
\end{equation*}
where the segments $s_{1:3}$ are pairwise different, and so are  the $q_{1:3}$ too.
Thanks to the Mecke-type formula for STIT tessellations, Theorem 3.1 of \cite{NNTW17},  there is an $\eta_M >0$, such that
\begin{multline*}
p_{\mathbf{ij}}(0) \leq  \frac{2\eta_M }{3! \cdot 3!}\int_{\mathcal{H}^6}\PPP{\partial Y_1\cap (B^o(H_{1:3})\cup B^o(H'_{1:3}))=\emptyset}\ind{(c(H_{1:3}), c(H'_{1:3}))\in \mathbf{i}\times \mathbf{j}}\ind{R(H_{1:3})> R(H'_{1:3})>v_\rho}\mathrm{d}H_{1:3}\mathrm{d}H'_{1:3}.
\end{multline*}
Here, the inequality is due to the omission of the indicators of the events $\{H_{1},H_2,H_3\}\cap B(q_{1:3})=\emptyset$ and $\{H'_{1},H'_2,H'_3\}\cap B(s_{1:3}) = \emptyset$. The factor 2 appears because we have added the indicator function $\ind{R(H_{1:3})> R(H'_{1:3})}$, which for a stationary and isotropic tessellation means no loss of generality. The factor $\eta_M$ comes from the application of the Mecke-type formula for the different possible arrangements of $s_{1:3}\in m_1^3$ and $q_{1:3}\in m_1^3$. We remark that the Mecke-type formula for maximal segments of a STIT tessellation is technically much more involved than the respective Mecke formula for a Poisson line tessellation. Here we omit a detailed derivation  of the formula above, as it would require quite a few pages; cf. the proofs in \cite{NNTW17}.

Now, let $\epsilon>0$ be fixed. For $H_{1:3}$ and $H'_{1:3}$ with $R(H_{1:3})> R(H'_{1:3})$ we write $B(c_1,r_1):=B(H_{1:3})$ and $B(c_2,r_2):=B(H'_{1:3})$. We consider two sub-cases. First, if $r_1\geq v_\rho(1+\epsilon)$, we use the fact that
 \[\PPP{\partial Y_1\cap (B(c_1,r_1)\cup B(c_2,r_2))=\emptyset}\leq \PPP{\partial Y_1\cap B(c_1,r_1)=\emptyset}.\] This yields
\[\PPP{\partial Y_1\cap (B(c_1,r_1)\cup B(c_2,r_2))=\emptyset}\leq e^{-2 r_1}.\]
Second, if $r_1<v_\rho(1+\epsilon)$, we use the inequality 
\[\PPP{\partial Y_1\cap (B(c_1,r_1)\cup B(c_2,r_2))=\emptyset}\leq \PPP{\partial Y_1\cap (B(c_1,r_2)\cup B(c_2,r_2))=\emptyset},\] which follows from the assumption that $r_1> r_2$.  Hence, according to Lemma \ref{Le:twoballs}, there exists a constant $\eta_1 >0$ such that 
\[\PPP{\partial Y_1\cap (B(c_1,r_1)\cup B(c_2,r_2))=\emptyset} \leq  \eta_1 \, e^{-2(1+\frac{2}{\pi}) r_2}.\] 
Summarizing, there exists a constant $\eta >0$ for which
\[p_{\mathbf{ij}}(0) \leq \eta \, (p_{\mathbf{ij}}^{> \epsilon}(0) + p_{\mathbf{ij}}^{<\epsilon}(0)),  \]
where
\[ p_{\mathbf{ij}}^{> \epsilon}(0):=  \int_{\mathcal{H}^6} e^{-2 R(H_{1:3})} \ind{(c(H_{1:3}), c(H'_{1:3}))\in \mathbf{i}\times \mathbf{j}}\ind{R(H_{1:3})> R(H'_{1:3})>v_\rho}\ind{R(H_{1:3})\geq v_\rho(1+\epsilon)}\mathrm{d}H_{1:3}\mathrm{d}H'_{1:3},\]
and 
\[ p_{\mathbf{ij}}^{< \epsilon}(0):= \int_{\mathcal{H}^6} e^{-2(1+\frac{2}{\pi}) R(H'_{1:3})}\ind{(c(H_{1:3}), c(H'_{1:3}))\in \mathbf{i}\times \mathbf{j}}\ind{v_\rho(1+\epsilon)>R(H_{1:3})> R(H'_{1:3})>v_\rho}\mathrm{d}H_{1:3}\mathrm{d}H'_{1:3}.\]
Now we provide  upper bounds for $p_{\mathbf{ij}}^{> \epsilon}(0)$ and $p_{\mathbf{ij}}^{< \epsilon}(0)$.

\paragraph{(i) Upper bound for $p_{\mathbf{ij}}^{> \epsilon}(0)$} We apply a Blaschke-Petkantschin type change of variables (see e.g. Theorem 7.3.2 of \cite{SW}). To do it we denote by  $a(u_{1:3})$ the area of the convex hull of $u_{1:3}\in (S^1)^3$. This gives
 \begin{align*}
 p_{\mathbf{ij}}^{> \epsilon}(0) & =  \int_{(S^1)^6} \int_{\mathbf{i}\times \mathbf{j}}\int_{\RR_+^2} a(u_{1:3})\ a(u'_{1:3})\ e^{-2 r_1}\ind{r_1> r_2>v_\rho}\ind{r_1\geq v_\rho(1+\epsilon)}   \mathrm{d}r_1\mathrm{d}r_2\mathrm{d}c_1\mathrm{d}c_2\sigma^{\otimes 3}(\mathrm{d}u_{1:3})\sigma^{\otimes 3}(\mathrm{d}u'_{1:3})\\
 & \leq \eta_2 \  a( \mathbf{i})^2 \ v_\rho\ e^{-2 v_\rho(1+\epsilon)},
 \end{align*}
for some positive constant $\eta_2$. This, together with (\ref{eq:defvrho}), and the fact that $a( \mathbf{i})=O(\log\log\rho)$, implies
\[p_{\mathbf{ij}}^{> \epsilon}(0)  = O\left( (\log\log\rho)^2\cdot \log\rho\cdot\rho^{-2(1+\epsilon)}  \right).\]

{\paragraph{(ii) Upper bound for $p_{\mathbf{ij}}^{< \epsilon}(0)$} Again we apply a Blaschke-Petkantschin type change of variables. This yields
 \begin{align*}
 p_{\mathbf{ij}}^{< \epsilon}(0) & =  \int_{(S^1)^6} \int_{\mathbf{i}\times \mathbf{j}}\int_{\RR_+^2} a(u_{1:3})\ a(u'_{1:3})\
e^{-2(1+\frac{2}{\pi}) r_2}\ind{ v_\rho(1+\epsilon)>r_1> r_2>v_\rho}   \mathrm{d}r_1\mathrm{d}r_2\mathrm{d}c_1\mathrm{d}c_2\sigma^{\otimes 3}(\mathrm{d}u_{1:3})\sigma^{\otimes 3}(\mathrm{d}u'_{1:3})\\
 & \leq \eta_3\  a( \mathbf{i})^2 \ v_\rho \ e^{-2(1+\frac{2}{\pi}) v_\rho},
 \end{align*}
for some positive constant $\eta_3$. In the same spirit as above, we get
\[p_{\mathbf{ij}}^{< \epsilon}(0)  = O\left( (\log\log\rho)^2\cdot \log\rho\cdot \rho^{-2(1+\frac{2}{\pi})}  \right).\]
Combining the upper bounds for $p_{\mathbf{ij}}^{> \epsilon}(0)$ and for $p_{\mathbf{ij}}^{< \epsilon}(0)$ we obtain for $0<\epsilon\leq \frac{2}{\pi}$ that
 \[p_{\mathbf{ij}}(0)= O\left( (\log\log\rho)^2\cdot \log\rho\cdot \rho^{-2(1+\epsilon)}  \right).\]

\subsubsection{Case $k=1$}
Now the calculation of an upper bound for $p_{\mathbf{ij}}(1)$  will be sketched only, because we proceed in the same spirit as in the case $k=0$. 
\begin{multline*}
p_{\mathbf{ij}}(1):=\\
 {\frac{1}{3! \cdot 2!}}\mathbb{E}\left[\sum_{s_{1:3}\in m_1^3}\ \sum_{q_{2:3}\in m_1^3 \setminus \{ s_{1:3}\} }\ind{R(s_{1:3})>v_\rho, R(s_1,q_{2:3})>v_\rho} 
\ind{\partial Y_1\cap (B^o(s_{1:3})\cup B^o(s_1,q_{2:3}))=\emptyset}\ind{(c(s_{1:3}), c(s_1,q_{2:3}))\in \mathbf{i}\times \mathbf{j}}\right].
\end{multline*}
Applying the Mecke-type formula for STIT tessellations, and discussing two sub-cases as above, namely $r_1\geq v_\rho(1+\epsilon)$ and $r_1 < v_\rho(1+\epsilon)$, we see that that there is an $\eta >0$ such that
\[p_{\mathbf{ij}}(1) \leq \eta \, (p_{\mathbf{ij}}^{> \epsilon}(1) + p_{\mathbf{ij}}^{<\epsilon}(1)),  \]
where
\[ p_{\mathbf{ij}}^{> \epsilon}(1):= \int_{\mathcal{H}^5} e^{-2 R(H_{1:3})} \ind{(c(H_{1:3}), c(H_1,H'_{2:3}))\in \mathbf{i}\times \mathbf{j}}\ind{ v_\rho< R(H_1,H'_{2:3})<     R(H_{1:3})}\ind{R(H_{1:3})\geq v_\rho(1+\epsilon)}\mathrm{d}H_{1:3}\mathrm{d}H'_{2:3},\]
and 
\[ p_{\mathbf{ij}}^{< \epsilon}(1):=  \int_{\mathcal{H}^5} e^{-2(1+\frac{2}{\pi}) R(H_1,H'_{2:3})}\ind{(c(H_{1:3}), c(H_1,H'_{2:3}))\in \mathbf{i}\times \mathbf{j}}\ind{v_\rho < R(H_1,H'_{2:3})  <R(H_{1:3}) <   v_\rho(1+\epsilon) }\mathrm{d}H_{1:3}\mathrm{d}H'_{2:3}.\]
Below, we provide   upper bounds for $p_{\mathbf{ij}}^{> \epsilon}(1)$ and $p_{\mathbf{ij}}^{< \epsilon}(1)$.

\paragraph{(i) Upper bound for $p_{\mathbf{ij}}^{> \epsilon}(1)$} Because $\ind{R(H_1,H'_{2:3})>v_\rho}\leq 1$, it follows from Fubini's theorem that 
\[ p_{\mathbf{ij}}^{> \epsilon}(1)\leq  \int_{\mathcal{H}^3} e^{-2R(H_{1:3})} \ind{c(H_{1:3})\in \mathbf{i}}\ind{R(H_{1:3})\geq v_\rho(1+\epsilon)}\left(  \int_{\mathcal{H}^2} \ind{c(H_1,H'_{2:3})\in \mathbf{j}} \ind{R(H_1,H'_{2:3})< R(H_{1:3})}\mathrm{d}H'_{2:3}\right)\mathrm{d}H_{1:3}.\]
Since $\mathbf{j}$ is a square with diameter $\sqrt{2}\frac{\sqrt{\pi \ \rho}}{ \left\lfloor \sqrt{\frac{\pi \ \rho}{\log\log\rho}} \right\rfloor }$, it follows from Lemma \ref{lem:G(H_1)}  that for all lines $H_1$, there is a constant ${\eta_4}>0$ such that
\[ \int_{\mathcal{H}^2} \ind{c(H_1,H'_{2:3})\in \mathbf{j}} \ind{R(H_1,H'_{2:3})< R(H_{1:3})}\mathrm{d}H'_{2:3} \leq {\eta_4} \  \log\log\rho.\]
Thus
\begin{align*}
p_{\mathbf{ij}}^{> \epsilon}(1) & \leq  {\eta_4} \ \log\log\rho \int_{\mathcal{H}^3} e^{-2 R(H_{1:3})} \ind{c(H_{1:3})\in \mathbf{i}}\ind{R(H_{1:3})\geq v_\rho(1+\epsilon)}\mathrm{d}H_{1:3}\\
& = {\eta_5} \, a(\mathbf{i})\,  \log\log\rho \int_{v_\rho(1+\epsilon)}^\infty e^{-2 r}\mathrm{d}r\\
& = O\left( (\log\log\rho)^{{2}}\cdot \rho^{-2(1+\epsilon)} \right),
\end{align*}
where the second line comes from the Blaschke-Petkantschin type change of variables, and ${\eta_5>0}$. 

\paragraph{(ii) Upper bound for $p_{\mathbf{ij}}^{< \epsilon}(1)$} This time we write
\begin{multline*} p_{\mathbf{ij}}^{< \epsilon}(1)\leq \int_{\mathcal{H}^3} e^{-2(1+\frac{2}{\pi}) R(H_1,H'_{2:3})} \ind{c(H_1,H'_{2:3})\in \mathbf{j}}\ind{v_\rho<R(H_1,H'_{2:3})}\\\times \left(  \int_{\mathcal{H}^2} \ind{c(H_{1:3})\in \mathbf{i}} \ind{R(H_{1:3})< v_\rho(1+\epsilon)}\mathrm{d}H_{2:3}\right)\mathrm{d}H_{1}\mathrm{d}H'_{2:3}.\end{multline*}
According to  Lemma \ref{lem:G(H_1)}, we have for some constant ${\eta_6}>0$
\[ \int_{\mathcal{H}^2} \ind{c(H_{1:3})\in \mathbf{i}} \ind{R(H_{1:3})< v_\rho(1+\epsilon)}\mathrm{d}H_{2:3}\leq {\eta_6}\  \log\log\rho.\]
Thus, for some constant ${\eta_7}>0$
\begin{align*}
p_{\mathbf{ij}}^{< \epsilon}(1)&\leq  {\eta_6} \, \log\log\rho \int_{\mathcal{H}^3} e^{-2(1+\frac{2}{\pi}) R(H_1,H'_{2:3})} \ind{c(H_1, H'_{2:3})\in \mathbf{j}}\ind{v_\rho<R(H_1,H'_{2:3})}\mathrm{d}H_1\mathrm{d}H'_{2:3}\\
& = {\eta_7} \,  a(\mathbf{j})\,  \log\log\rho \int_{v_\rho}^\infty e^{-2(1+\frac{2}{\pi}) r}dr\\
& = O\left( (\log\log\rho)^{{2}} \cdot  \rho^{-2(1+\frac{2}{\pi})} \right),   
\end{align*}
where the second line again comes from the Blaschke-Petkantschin type change of variables.  Summarizing,  we obtain that $p_{\mathbf{ij}}(1)= O\left( {(\log\log\rho)^{{2}}} \cdot  \rho^{-2(1+\epsilon)} \right)$ for $0<\epsilon \leq \frac{2}{\pi}$.

\subsubsection{Case $k=2$} 
Proceeding now exactly along the same lines as above, we see that there is an $\eta>0$ such that
\[p_{\mathbf{ij}}(2) \leq \eta \, (p_{\mathbf{ij}}^{> \epsilon}(2) + p_{\mathbf{ij}}^{<\epsilon}(2)),  \]
where
\[ p_{\mathbf{ij}}^{> \epsilon}(2):= \int_{\mathcal{H}^4} e^{-2 R(H_{1:3})} \ind{(c(H_{1:3}), c(H_{1:2},H'_{3}))\in \mathbf{i}\times \mathbf{j}}\ind{R(H_{1:3})> R(H_{1:2},H'_{3})>v_\rho}\ind{R(H_{1:3})\geq v_\rho(1+\epsilon)}\mathrm{d}H_{1:3}\mathrm{d}H'_{3},\]
and 
\[ p_{\mathbf{ij}}^{< \epsilon}(2):= \int_{\mathcal{H}^4} e^{-2(1+\frac{2}{\pi}) R(H_{1:2},H'_{3})}\ind{(c(H_{1:3}), c(H_{1:2},H'_{3}))\in \mathbf{i}\times \mathbf{j}}\ind{ v_\rho(1+\epsilon)> R(H_{1:3})>R(H_{1:2},H'_{3})>v_\rho }\mathrm{d}H_{1:3}\mathrm{d}H'_{3}.\]

\paragraph{(i) Upper bound for $p_{\mathbf{ij}}^{> \epsilon}(2)$} We write
\[ p_{\mathbf{ij}}^{> \epsilon}(2)\leq  \int_{\mathcal{H}^3} e^{-2 R(H_{1:3})} \ind{c(H_{1:3})\in \mathbf{i}}\ind{R(H_{1:3})\geq v_\rho(1+\epsilon)}\left(  \int_{\mathcal{H}} \ind{c(H_{1:2},H'_{3})\in \mathbf{j}} \ind{R(H_{1:2},H'_{3})< R(H_{1:3})}\mathrm{d}H'_3\right)\mathrm{d}H_{1:3}.\]
Since $\mathbf{j}$ is a square with diameter $\sqrt{2}\frac{\sqrt{\pi\ \rho}}{ \left\lfloor \sqrt{\frac{\pi\ \rho}{\log\log\rho}} \right\rfloor }$, it follows from Lemma \ref{lem:G(H_2)} that for all $H_{1:3}$, there is a constant ${\eta_8} >0$ such that
\[ \int_{\mathcal{H}} \ind{c(H_{1:2},H'_{3})\in \mathbf{j}} \ind{R(H_{1:2},H'_{3})< R(H_{1:3})}\mathrm{d}H'_{3} \leq {\eta_8} \, \sqrt{\log\log\rho }.\]
Proceeding along the same lines as above, we obtain
\begin{align*}
p_{\mathbf{ij}}^{> \epsilon}(2) = O\left(\left(\log\log\rho \right)^{{\frac{3}{2}}}\cdot  \rho^{-2(1+\epsilon)}  \right).
\end{align*}

\paragraph{(ii) Upper bound for $p_{\mathbf{ij}}^{< \epsilon}(2)$} Now
\[ p_{\mathbf{ij}}^{< \epsilon}(2)\leq  \int_{\mathcal{H}^3} e^{-2(1+\frac{2}{\pi})R(H_{1:2},H'_{3})} \ind{c(H_{1:2},H'_{3})\in \mathbf{j}}\ind{v_\rho<R(H_{1:2},H'_{3})}\left(  \int_{\mathcal{H}} \ind{c(H_{1:3})\in \mathbf{i}} \ind{R(H_{1:3})< v_\rho(1+\epsilon)}\mathrm{d}H_{3}\right)\mathrm{d}H_{1:2}\mathrm{d}H'_3 .\]
This yields
\begin{align*}
p_{\mathbf{ij}}^{< \epsilon}(2) = O\left(\left(\log\log\rho \right)^{{\frac{3}{2}}}\cdot  \rho^{-2(1+\frac{2}{\pi})}  \right).
\end{align*}
Summarizing, we obtain $p_{\mathbf{ij}}(2)= O\left( {(\log\log\rho)^{{\frac{3}{2}}}} \cdot  \rho^{-2(1+\epsilon)} \right)$ for $0<\epsilon \leq \frac{2}{\pi}$.

\medskip

Now, having upper bounds for  $p_{\mathbf{ij}}(k)$, $k=0,1,2$, we can complete the proof that $b_2$ converges to 0 as $\rho \to \infty$. Indeed,  {let} $\epsilon \in (0,\frac{2}{\pi}]$ {and} $\beta \in (0,1)$ {be fixed}.  {Let}  $\epsilon' :=2(1+\epsilon) - (1+\beta)$. Then
 \[b_2 =\sum_{\mathbf{i}\in V}\ \sum_{\mathbf{i}\neq \mathbf{j}\in S(\mathbf{i},\rho^{\beta /2} )}p_{\mathbf{ij}}
\leq |V|\cdot |S(\mathbf{i},\rho^{\beta /2} )|\cdot \sum_{k=0}^2 p_{\mathbf{ij}} (k)
= O\left( {\log\log\rho\cdot \log\rho\cdot \rho^{-\epsilon'}}\right).\]
{Since $\epsilon'>0$, this proves that $b_2$ converges to 0 as $\rho$ goes to infinity}.

\subsection{Proof of the assertion that $b_3$ converges to 0 for $\rho \to \infty$}
\label{subsection:mixing}
First, we rewrite $b_3$ as a series. To do it, let $\rho>\rho_0(\tau)$ and let $\mathbf{i}\in V$. According to \eqref{eq:rhotau}, we know that  $\sum_{\mathbf{j}\not\in S(\mathbf{i},\rho^{\beta /2} )} \ind{M_\mathbf{j}>v_\rho}=N_{C(\mathbf{i})}(v_\rho)$, where $C(\mathbf{i})$ is defined in \eqref{eq:defCi}. Now, using the Bayes theorem,  elementary calculations on conditional probabilities yield
\begin{align*}
\label{eq:majbprime3}
b_3 & = \sum_{\mathbf{i}\in V}p_\mathbf{i}\sum_{k=0}^\infty \left| \PPP{N_{C(\mathbf{i})}(v_\rho)=k| M_\mathbf{i}>v_\rho } - \PPP{N_{C(\mathbf{i})}(v_\rho)=k}  \right| \\
       & \leq   \left( \sum_{\mathbf{i}\in V}p_\mathbf{i}  \right)\cdot \left( \sup_{\mathbf{i}\in V}\sum_{k=0}^\infty \left| \PPP{N_{C(\mathbf{i})}(v_\rho)=k| M_\mathbf{i} >v_\rho} - \PPP{N_{C(\mathbf{i})}(v_\rho)=k}  \right| \right).
\end{align*}
Moreover, according to \eqref{Campbellmark} and \eqref{eq:pi}, we have $\sum_{\mathbf{i}\in V}p_\mathbf{i}= \tau$. Thus, to prove that $b_3$ converges to 0, it is sufficient to prove that $b'_3$ converges to 0, where
\begin{equation*}
b'_3:=\sup_{\mathbf{i}\in V}\sum_{k=0}^\infty \left| \PPP{N_{C(\mathbf{i})}(v_\rho)=k| M_\mathbf{i} >v_\rho} - \PPP{N_{C(\mathbf{i})}(v_\rho)=k}  \right|.
\end{equation*}

The main idea is to apply Lemmas \ref{le:encapsulation} and \ref{le:measurable}. To be in the framework of these lemmas, we intersect the events appearing in the series  with the events $ \{M_{C(\mathbf{i})}\leq v_\rho^2\}$ and $ \{M_\mathbf{i}\leq v_\rho^2\}$ respectively. Notice that these two events occur with high probability. With standard computations, we can prove that for each 4-tuple of events $A,B,C,D$, with $\PPP{B\cap D}\neq 0$, we have 
\[|\PPP{A|B} - \PPP{A}| \leq \left| \PPP{A\cap C|B\cap D} - \PPP{A\cap C} \right| + \frac{1}{\PPP{B}}(2\PPP{A\cap C^c}+\PPP{A\cap D^c}+2\PPP{A}\PPP{D^c}).\]

Applying the above inequality to the events $A=\{N_{C(\mathbf{i})}(v_\rho)=k\}$, $B=\{M_\mathbf{i}>v_\rho\}$, $C=\{M_{C(\mathbf{i})}\leq v_\rho^2\}$ and $D=\{M_\mathbf{i}\leq v_\rho^2\}$, we obtain  for all  $\mathbf{i}$ and all $k\in \NN_0$: 
\begin{equation*}
 \left| \PPP{N_{C(\mathbf{i})}(v_\rho)=k| M_\mathbf{i} >v_\rho} - \PPP{N_{C(\mathbf{i})}(v_\rho)=k}  \right| \leq p_{k,\mathbf{i}} + q_{k,\mathbf{i}},
 \end{equation*}
 where
 \begin{equation*}
  p_{k,\mathbf{i}}=\left| \PPP{\left.\{N_{C(\mathbf{i})}(v_\rho)=k\}\cap \{M_{C(\mathbf{i})}\leq v_\rho^2\}\right| v_\rho < M_\mathbf{i} \leq v_\rho^2} - \PPP{\{N_{C(\mathbf{i})}(v_\rho)=k\}\cap \{M_{C(\mathbf{i})}\leq v_\rho^2\}}  \right|,
 \end{equation*}
and
\begin{multline*}
 q_{k,\mathbf{i}} = \frac{1}{\PPP{M_\mathbf{i}>v_\rho}} \Big(2\PPP{\{N_{C(\mathbf{i})}(v_\rho)=k  \}\cap \{M_{C(\mathbf{i})}>v_\rho^2\}} + \PPP{\{N_{C(\mathbf{i})}(v_\rho)=k\}\cap \{M_\mathbf{i}>v_\rho^2\}}\\
  + 2\PPP{N_{C(\mathbf{i})}(v_\rho)=k}\PPP{M_\mathbf{i}>v_\rho^2} \Big).
\end{multline*}
Thus, to prove that $b'_3$ converges to 0, we have to prove that $\sup_{\mathbf{i}\in V}\sum_{k=0}^\infty p_{k, \mathbf{i}}$ and $\sup_{\mathbf{i}\in V}\sum_{k=0}^\infty q_{k, \mathbf{i}}$ converge to 0.  This is done in Subsections \ref{subsec:pk} and \ref{subsec:qk}. We begin with $\sup_{\mathbf{i}\in V}\sum_{k=0}^\infty q_{k, \mathbf{i}}$ because it is the simplest case.

\subsubsection{Proof of the assertion that $\sup_{\mathbf{i}\in V}\sum_{k=0}^\infty q_{k, \mathbf{i}}$ converges to 0}
\label{subsec:qk}
Roughly, the fact that $\sup_{\mathbf{i}\in V}\sum_{k=0}^\infty q_{k, \mathbf{i}}$ converges to 0  comes from the fact that the events $ \{M_{C(\mathbf{i})}\leq v_\rho^2\}$ and $ \{M_\mathbf{i}\leq v_\rho^2\}$  occur with high probability. Let $\mathbf{i}\in V$. Summing over $k\in \NN_0$, we have
\[\sum_{k=0}^\infty q_{k, \mathbf{i}} = \frac{1}{\PPP{M_\mathbf{i}>v_\rho}} \left(2\PPP{M_{C(\mathbf{i})}>v_\rho^2} + 3\PPP{M_\mathbf{i}>v_\rho^2} \right).\]
Now, analogously to \eqref{eq:pi}, we know that
\begin{equation*}
\PPP{M_{C(\mathbf{i})}>v_\rho^2} \leq \EEE {\sum_{\mathbf{j}\in C(\mathbf{i})} \ind{M_\mathbf{j}> v_\rho^2}} = \gamma_1 a(C(\mathbf{i}))\ \PPP{R(\mathcal{Z})>v_\rho^2}.
\end{equation*} 
Furthermore, the stationarity of the tessellation yields that $\PPP{M_\mathbf{i}>v_\rho^2}\leq \PPP{M_{C(\mathbf{i})}>v_\rho^2}$ for all  $\rho$ for which $a(\mathbf{i})\leq a(C(\mathbf{i}))$. {Thus   
\begin{equation*}
\sum_{k=0}^\infty q_{k, \mathbf{i}} 
\leq   \frac{5\PPP{M_{C(\mathbf{i})}>v_\rho^2}}{\PPP{M_\mathbf{i}>v_\rho}} 
\leq \frac{5 a(C(\mathbf{i}))\PPP{R(\mathcal{Z})>v_\rho^2}}{a(\mathbf{i})\PPP{R(\mathcal{Z})>v_\rho}}.
\end{equation*}
Since $a(C(\mathbf{i})) \leq a(W_\rho )$ and $a(\mathbf{i})=\frac{a(W_\rho)}{|V|}$, we obtain from \eqref{eq:subdivision}  that
\[\sum_{k=0}^\infty q_{k, \mathbf{i}} \leq \frac{5\pi\rho}{\log\log \rho}\, \frac{\PPP{R(\mathcal{Z})>v_\rho^2}}{\PPP{R(\mathcal{Z})>v_\rho}}.\]
}
An application of \eqref{eq:exponential} and  \eqref{eq:defvrho} yields
\[\sum_{k=0}^\infty q_{k, \mathbf{i}} = O\left(  (\log\log \rho )^{-1} \, \rho^{-\frac{1}{2}\log \rho +\log \tau +2} \right).\]}

As the last term does not depend on $\mathbf{i}\in V$, this is also an upper bound for $\sup_{\mathbf{i}\in V}\sum_{k=0}^\infty q_{k, \mathbf{i}}$. Obviously it converges to 0 as $\rho \to \infty$.

\subsubsection{Proof of the assertion that $\sup_{\mathbf{i}\in V}\sum_{k=0}^\infty p_{k, \mathbf{i}}$ converges to 0}
\label{subsec:pk}
Let $\mathbf{i}\in V$ be fixed. For each $0<s< 1$,  we denote by $M_\mathbf{i}^{[1-s]}$ the maximum  of the inradii over all cells with nucleus in $\mathbf{i}$ for a STIT tessellation at time $1-s$. When the time is not specified, the underlying STIT tessellation is at time $1$, e.g. $M_\mathbf{i}$ denotes the same maximum but this time for a STIT tessellation at time $1$, in accordance with \eqref{eq:notationmaxstit}. Now taking 
 \[K=\left(\bigcup_{\mathbf{j}\in S(\mathbf{i},\rho^{\beta /2} )} \mathbf{j}\right) \ominus \mathbf{S}_0 \quad \text{and} \quad K'=\mathbf{i}\oplus \mathbf{S}_0,\] 
it follows from Lemmas \ref{le:encapsulation} and \ref{le:measurable} that for each $s\in (0,1)$, 
\begin{multline*}
p_{k,\mathbf{i}}\leq  \PPP{\{ N_{C(\mathbf{i})}(v_\rho) = k\}\cap \{M_{C(\mathbf{i})}\leq v_\rho^2\}| v_\rho <M_\mathbf{i}\leq v_\rho^2}  \\ 
\begin{split}
&- \frac{\PPP{ v_\rho <{M_\mathbf{i}^{[1-s]}}\leq v_\rho^2}}{\PPP{v_\rho <M_\mathbf{i}\leq v_\rho^2}}\cdot \PPP{ N_{C(\mathbf{i})}(v_\rho) = k, M_{C(\mathbf{i})}\leq v_\rho^2, S(K,K')<s, Y_s\wedge K'=K'}\\
& + \left| \frac{\PPP{ v_\rho <{M_\mathbf{i}^{[1-s]}}\leq v_\rho^2}}{\PPP{ v_\rho <M_\mathbf{i}\leq v_\rho^2}} - 1  \right|\cdot \PPP{ N_{C(\mathbf{i})}(v_\rho) = k, M_{C(\mathbf{i})}\leq v_\rho^2}\\
&  + \frac{\PPP{v_\rho <{M_\mathbf{i}^{[1-s]}}\leq v_\rho^2}}{\PPP{ v_\rho <M_\mathbf{i}\leq v_\rho^2}}\cdot \PPP{ \{ N_{C(\mathbf{i})}(v_\rho) = k, M_{C(\mathbf{i})}\leq v_\rho^2\}\cap \{S(K,K')<s, Y_s\wedge K'=K'\}^c}.
\end{split}
\end{multline*}
Notice that in the above equation we have considered the sets $\mathbf{i}$ and $C(\mathbf{i})$ instead of the sets $\mathbf{i}^\circ$ and $C(\mathbf{i})^\circ$, which appear in Lemma \ref{le:measurable}. Actually, this does not modify the probability which are considered because a.s. no incenter of a cell is located at the boundary of ${\mathbf{i}}$ or of $C(\mathbf{i})$ respectively. Now, summing over $k\geq 0$ and using the fact that $\sum_{k=0}^\infty \PPP{N_{C(\mathbf{i})}(v_\rho) = k} = 1 $, we have 
\begin{multline*}
\sum_{k=0}^\infty p_{k, \mathbf{i}}\leq  \PPP{M_{C(\mathbf{i})}\leq v_\rho^2|v_\rho <M_\mathbf{i}\leq v_\rho^2} \\
\begin{split} 
& -  \frac{\PPP{ v_\rho <{M_\mathbf{i}^{[1-s]}} \leq v_\rho^2}}{\PPP{v_\rho <M_\mathbf{i}\leq v_\rho^2}}\cdot \PPP{M_{C(\mathbf{i})}\leq v_\rho^2, S(K,K')<s, Y_s\wedge K'=K'}\\
& + \left| \frac{\PPP{ v_\rho <{M_\mathbf{i}^{[1-s]}}\leq v_\rho^2}}{\PPP{v_\rho <M_\mathbf{i}\leq v_\rho^2}} - 1  \right|\cdot \PPP{M_{C(\mathbf{i})}\leq v_\rho^2} \\ 
&  + \frac{\PPP{ v_\rho <{M_\mathbf{i}^{[1-s]}}\leq v_\rho^2}}{\PPP{v_\rho <M_\mathbf{i}\leq v_\rho^2}}\cdot \PPP{   \{M_{C(\mathbf{i})}\leq v_\rho^2\}\cap \{S(K,K')<s, Y_s\wedge K'=K'\}^c}.
\end{split}
\end{multline*}
{Bounding the probabilities} $ \PPP{M_{C(\mathbf{i})}\leq v_\rho^2|v_\rho <M_\mathbf{i}\leq v_\rho^2}$ and  $\PPP{M_{C(\mathbf{i})}\leq v_\rho^2}$ {by 1 respectively, and} using the trivial inequalities
\begin{equation*} \PPP{   \{M_{C(\mathbf{i})}\leq v_\rho^2\}\cap \{S(K,K')<s, Y_s\wedge K'=K'\}^c}\leq  \PPP{ \{S(K,K')<s, Y_s\wedge K'=K'\}^c}
\end{equation*}
and 
\begin{equation*} \PPP{M_{C(\mathbf{i})}\leq v_\rho^2, S(K,K')<s, Y_s\wedge K'=K'}\geq \PPP{S(K,K')<s, Y_s\wedge K'=K'} - \PPP{    M_{C(\mathbf{i})}> v_\rho^2},
\end{equation*}
we obtain 
\begin{multline*}
\sum_{k=0}^\infty p_{k, \mathbf{i}} \leq   1- \frac{\PPP{v_\rho <{M_\mathbf{i}^{[1-s]}}\leq v_\rho^2}}{\PPP{v_\rho <M_\mathbf{i}\leq v_\rho^2}}\cdot \left( \PPP{S(K,K')<s, Y_s\wedge K'=K'} - \PPP{    M_{C(\mathbf{i})}> v_\rho^2}
\right)\\
\begin{split}
& + \left| \frac{\PPP{v_\rho <{M_\mathbf{i}^{[1-s]}}\leq v_\rho^2}}{\PPP{v_\rho <M_\mathbf{i}\leq v_\rho^2}} - 1  \right|   + \frac{\PPP{ v_\rho <{M_\mathbf{i}^{[1-s]}}\leq v_\rho^2}}{\PPP{v_\rho <M_\mathbf{i}\leq v_\rho^2}}\cdot \PPP{  \{S(K,K')<s, Y_s\wedge K'=K'\}^c}.
\end{split}
\end{multline*}
Analogously to \eqref{eq:pi}, and applying  \eqref{eq:exponential}, \eqref{eq:defvrho}, and the fact that $a(C(\mathbf{i}))=O\left( \rho  \right)$ we obtain that
\begin{equation*}
\PPP{M_{C(\mathbf{i})}>v_\rho^2}  \leq   \EEE {\sum_{\mathbf{j}\in C(\mathbf{i})} \ind{M_\mathbf{j}> v_\rho^2}} = \gamma_1\,  a(C(\mathbf{i}))\ \PPP{R(\mathcal{Z})>v_\rho^2} 
 = O(\rho^{-\frac{1}{2}\log \rho +\log \tau +1}),
\end{equation*} 
which converges to 0 as $\rho$ goes to infinity. Thus, to prove that $\sum_{k=0}^\infty p_{k, \mathbf{i}}$ converges to 0, it is sufficient to choose $s$, as a function of $\rho$, in such a way that the following properties hold: 
\begin{equation}\label{eq:aim1} \frac{\PPP{ v_\rho <{M_\mathbf{i}^{[{1}-s]}}\leq v_\rho^2}}{\PPP{v_\rho <M_\mathbf{i}\leq v_\rho^2}} \conv[\rho]{\infty}1\end{equation}
and 
\begin{equation}  \label{eq:aim2} \PPP{S(K,K')<s, Y_s\wedge K'=K'}\conv[\rho]{\infty}1.\end{equation}
To do it, we choose $s=\rho^{-\delta}$ with $0<\delta< \beta /2$.

First, we prove \eqref{eq:aim1}. Analogously to \eqref{eq:pi}, for $u=1,1-s$ and $k=1,2$, we have $\PPP{{M_\mathbf{i}^{[u]}}>v_\rho^k} = a(\mathbf{i})\gamma_{u} e^{-2uv_\rho^k}$, with $\rho > \rho_0(\tau)$. Since $\gamma_u=\frac{u^2}{\pi}$, this gives
\[ \frac{\PPP{ v_\rho <{M_\mathbf{i}^{[1-s]}}\leq v_\rho^2}}{\PPP{v_\rho <M_\mathbf{i}\leq v_\rho^2}}
=\frac{[1-s]^2 \left(  e^{-2 [1-s]v_\rho} - e^{-2 [1-s]v^2_\rho} \right)}{e^{-2 v_\rho} - e^{-2 v^2_\rho}}.\] 
This together with the fact that  $s\, v^2_\rho \conv[\rho]{\infty}0$ shows \eqref{eq:aim1}. 

\begin{center}
\begin{figure}
\begin{center}
\includegraphics[scale=0.5]{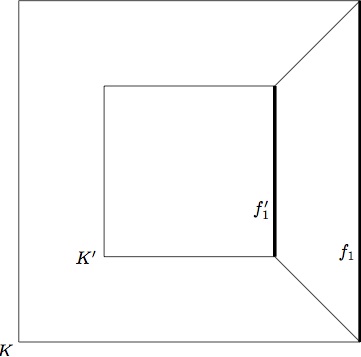}
\end{center}
\caption{The sides $f_1$, $f'_1$ of the squares $K$, $K'$, and the associated trapezoid}
\label{fig:separating}
\end{figure}
\end{center}

It remains to prove \eqref{eq:aim2}. To do it, we apply an inequality established in \cite{MN}. We first recall the framework of this paper. Let {$a$} be the length of the side of the square (centered at the origin) $K=\left(\bigcup_{\mathbf{j}\in S(\mathbf{i},\rho^{\beta /2} )} \mathbf{j}\right) \ominus \mathbf{S}_0$. It is clear that \begin{equation}\label{eq:estimatea} a=\sqrt{a(\mathbf{i})}(2\lfloor \rho^{\frac{\beta}{2}}\rfloor+1)-2v_\rho^2
 \end{equation}for some constant $c$. Let  $f_1= [-\frac{a}{2},\frac{a}{2}]\times \{\frac{a}{2}\}$, $f_2=\{\frac{a}{2}\}\times [-\frac{a}{2},\frac{a}{2}]$, $f_3=-f_1$ and $f_4=-f_2$ be the sides of $K$. Similarly, with the same orientation, we denote by $f'_1,f'_2,f'_3,f'_3$ the sides of the square  (centered at the origin) $K'=\mathbf{i}\oplus \mathbf{S}_0$. The length of each side of $K'$ is 
\begin{equation} \label{eq:estimatea'}   a'=\sqrt{a(\mathbf{i})}+2v_\rho^2. 
\end{equation} Now, let $L(K,K'):=\min\{\Lambda ([f_i|f'_i]): i=1,\ldots, 4\}$, as in equation (13) of \cite{MN}. Notice that the quantities $\Lambda ([f_i|f'_i])$ appearing in this minimum actually do not depend on $i$ because $\Lambda$ is invariant under rotation. Equation (18) of \cite{MN} states the following inequality:
\[ \PPP{S(K,K')<s, Y_s\wedge K'=K'} \geq e^{-s\Lambda[K']}\left( 1-e^{-sL(K,K')}\right)^4.\] 
To calculate $L(K,K')$, consider the trapezoid spanned by $f_1$ and $f'_1$ (see Figure \ref{fig:separating}). A line $H$ separates $f_1$ and $f'_1$, if and only if it intersects the two nonparallel sides of the trapezoid. Denote by $\ell$ the length of a diagonal of the trapezoid. The distance between $f_1$ and $f'_1$ is $\frac{1}{2}(a-a')$. A formula for the measure of all lines intersecting two planar compact convex sets can be found in \cite{Sa}, p. 33. Applying this to the linear segments $f_1$, $f'_1$, we obtain with some elementary geometry that \label{page:LKK}
\begin{equation*}
L(K,K')  =2 \ell -(a+a')=\sqrt{2}\, \sqrt{a^2 + a'^2} - (a+a') \geq (\sqrt{2}-1)a-a'.
\end{equation*}
According to \eqref{eq:estimatea} and \eqref{eq:estimatea'}, this gives {
\begin{align*}
L(K,K') & \geq (\sqrt{2}-1)\left( \sqrt{a(\mathbf{i})} (2\lfloor \rho^{\frac{\beta}{2}}\rfloor +1)-2v_\rho^2  \right) - (\sqrt{a(\mathbf{i})}+2v_\rho^2)  \\
& {= \left((\sqrt{2}-1)  (2\lfloor \rho^{\frac{\beta}{2}}\rfloor +1) -1 \right)\sqrt{a(\mathbf{i})} -2\sqrt{2}v_\rho^2}\\
& {\geq 2(\sqrt{2}-1)\rho^{\frac{\beta}{2}}-\sqrt{2}-2\sqrt{2}v_\rho^2,}
\end{align*}
where the second inequality comes from the facts that $\sqrt{a(\mathbf{i})}\geq 1$ (when $\rho$ is sufficiently large) and $\lfloor \rho^{\frac{\beta}{2}}\rfloor \geq \rho^{\frac{\beta}{2}}-1$. 
Using the fact that $\sqrt{2}\leq (\sqrt{2}-1)\rho^{\frac{\beta}{2}}$ (for $\rho$ sufficiently large), we get
\begin{align*}
L(K,K') & \geq 2(\sqrt{2}-1)\rho^{\frac{\beta}{2}}-(\sqrt{2}-1)\rho^{\frac{\beta}{2}}-2\sqrt{2}v_\rho^2\\
& = (\sqrt{2}-1)\rho^{\frac{\beta}{2}}-2\sqrt{2}v_\rho^2.
\end{align*}
}
Since $s=\rho^{-\delta}$ with $0<\delta< \beta /2$, we obtain that $s\, L(K,K') \conv[\rho]{\infty} \infty$.
 Furthermore, $s\Lambda[K'] = \frac{4}{\pi} s\, a'$, which converges to 0 as $\rho$ goes to infinity.  This shows \eqref{eq:aim2}. Consequently,  $b_3$ converges to 0 as $\rho$ goes to infinity, which completes the proof of Theorem \ref{Th:extremes}.

\section{Concluding remarks}
\label{sec:conclusion}
In this section, we discuss some possible extensions of Theorem \ref{Th:extremes}. First, in the proof of our main theorem, it is essential that the considered random tessellation  is stationary. Such a condition is standard in Stochastic Geometry, and is classical in Extreme Value Theory. Another assumption is that our STIT tessellation is isotropic. This is used in the proof of Lemma \ref{Le:twoballs} and in the computation of $L(K,K')$ (see p. \pageref{page:LKK}). However, our results can be extended to non-isotropic STIT tessellations, as long as the directional distribution of the dividing lines are atomless, and hence no parallel sides occur in the tessellation. 
 The asymptotic behaviour of the largest order statistics is the same as in Theorem 1.1, especially because \eqref{eq:exponential} remains true for stationary but non-isotropic STIT tessellations. If there are cells with pairs of parallel sides, their incircles must not be unique, and a different approach to the problem is required.

Our main theorem is given specifically for the two dimensional case with a fixed square-shaped window $W_\rho$ in order to keep our calculations simple. However, Theorem 1.1 remains true when the window is any convex body with non-empty interior. We conjecture that our results concerning the largest order statistics may be extended to higher dimensions. 

We now compare our method to the one used in \cite{ChenHem}. As mentioned on p. \pageref{page:moment}, our result is similar to Theorem 1.1 (ii) of \cite{ChenHem}. However, the proof of this theorem is based on the method of moments, which leads to very technical  computations of combinatorics, and does not provide rate of convergence for the Poisson approximation. On the opposite, our proof based on the Chen-Stein method, is more accurate in the sense that the rate of convergence can be made explicit. In Theorem 1.1 (i) of \cite{ChenHem}, the smallest order statistics for the inradius are also considered. We did not investigate this problem in the context of STIT tessellations because it relies on a simple adaptation of \cite{ChenHem}. The method is similar and is mainly based on the study of $U$-statistics, which are investigated  in \cite{ST}.

\paragraph{Acknowledgements} This work was partially supported by the French ANR grant ASPAG  (ANR-17-CE40-0017).


\begin{thebibliography}{99}

\bibitem{AGG} R. Arratia, L. Goldstein, and L. Gordon. Two moments suffice for {P}oisson approximations: the {C}hen-{S}tein method, \textit{Ann. Probab.,} (17) 9--25, 1989.

\bibitem{CC} P. Calka and N. Chenavier.  Extreme values for characteristic radii of a  {P}oisson-{V}oronoi tessellation, \textit{Extremes,} (3) 359--385, 2014.

\bibitem{Chen} N. Chenavier. A general study of extremes of stationary tessellations with examples. \textit{Stochastic Process. Appl.}, 124(9):2917–-2953, 2014.

\bibitem{ChenHem} N. Chenavier and R. Hemsley. Extremes for the inradius in the Poisson line tessellation. \textit{Adv. in Appl. Probab.}, 48(2):544–-573, 2016.

\bibitem{HF} L. de Haan and A. Ferreira. Extreme value theory. An introduction. \textit{Springer Series in Operations Research and Financial Engineering}. Springer, New York, 2006.

\bibitem{lachiezerey} R. Lachi\`eze-Rey. Mixing properties for STIT tessellations. \textit{Adv. in Appl. Probab.}, 43(1):40-–48, 2011.

\bibitem{LP} V. Lazarus and L. Pauchard. From craquelures to spiral crack patterns: influence of layer thickness on the crack patterns induced by desiccation. \textit{Soft Matter}, 7:2552-–2559, 2011.

\bibitem{MN3} S. Martinez and W. Nagel. Ergodic description of STIT tessellations. \textit{Stochastics}, 84(1):113–-134, 2012.

\bibitem{MN2} S. Martinez and W. Nagel. STIT tessellations have trivial tail $\sigma$-algebra. \textit{Adv. in Appl. Probab.}, 46(3):643–-660, 2014.

\bibitem{MN} S. Martinez and W. Nagel. The $\beta$-mixing rate of STIT tessellations. \textit{Stochastics }88, 3:396-–414, 2016.

\bibitem{mnw08a} J. Mecke, W. Nagel, and V. Weiss. A global construction of homogeneous random planar tessellations that are stable under iteration. \textit{Stochastics}, 80(1):51–-67, 2008.

\bibitem{NagelWeiss11} J. Mecke, W. Nagel, and V. Weiss. Some distributions for I-segments of planar random homogeneous STIT tessellations. \textit{Math. Nachr.}, 284(11-12):1483–-1495, 2011.

\bibitem{MoMa} L. Mosser and S. Matthai. Tessellations stable under iteration: Evaluation of application as an improved stochastic discrete fracture modeling algorithm. \textit{International Discrete Fracture Network Engineering Conference}, 2014.

\bibitem{NNTW17} W. Nagel, N. L. Nguyen, C. Th\"ale, and V. Weiss. A Mecke-type formula and Markov properties for STIT tessellation processes. \textit{ALEA Lat. Am. J. Probab. Math. Stat.}, 14(2):691-–718, 2017.

\bibitem{NagelWeiss05} W. Nagel and V. Weiss. Crack STIT tessellations: characterization of stationary random tessellations stable with respect to iteration. \textit{Adv. in Appl. Probab.}, 37(4):859–-883, 2005.


\bibitem{NagelWeiss08} W. Nagel and V. Weiss. Mean values for homogeneous STIT tessellations in 3d. \textit{Image Anal. Stereol.}, 27(1):29-–37, 2008.

\bibitem{Sa} L. A. Santal\'o. Integral geometry and geometric probability. Cambridge Mathematical Library. \textit{Cambridge University Press, Cambridge}, second edition, 2004. With a foreword by Mark Kac.

\bibitem{SW} R. Schneider and W. Weil. Stochastic and integral geometry. Probability and its Applications (New York). \textit{Springer-Verlag}, Berlin, 2008.

\bibitem{ST} M. Schulte and C. Th\"ale. The scaling limit of Poisson-driven order statistics with applications in geo- metric probability. \textit{Stochastic Process. Appl.}, 122(12):4096–-4120, 2012.

\bibitem{SKM} D. Stoyan, W. Kendall, and J. Mecke. Stochastic Geometry and its applications. \textit{Wiley}, 2008.


  

\end{thebibliography}

\end{document}